\documentclass[12pt,dvips]{amsart}
\usepackage{amsfonts, amssymb, latexsym, epsfig,
amsmath, amscd}
\usepackage{url}
\usepackage{tikz}
\usepackage{graphicx}
\usepackage{enumerate}
\usepackage[all]{xy}

\newtheorem{Theorem}{Theorem}[section]
\newtheorem*{theorem*}{Theorem}
\newtheorem{Proposition}[Theorem]{Proposition}
\newtheorem{Lemma}[Theorem]{Lemma}
\newtheorem{Claim}[Theorem]{Claim}
\newtheorem{Subclaim}[Theorem]{Subclaim}

\newtheorem{Algorithm}[Theorem]{Algorithm}
\newtheorem{Corollary}[Theorem]{Corollary}
\newtheorem{Problem}[Theorem]{Problem}
\newtheorem{Definition-Proposition}[Theorem]{Definition-Theorem}
\newtheorem{Main Conjecture}[Theorem]{Main Conjecture}
\newtheorem{Observation}[Theorem]{Observation}
\newtheorem{Conjecture}[Theorem]{Conjecture}
\newtheorem{Definition}[Theorem]{Definition}
\theoremstyle{remark}
\newtheorem{Example}[Theorem]{Example}
\newtheorem{prop}{Proposition}[section]
\newtheorem{Remark}[prop]{Remark}

\newcommand\union{\bigcup}

\newcommand{\caA}{\mathcal{A}}

\newcommand{\caI}{\mathcal{I}}

\newcommand{\caN}{\mathcal{N}}
\newcommand{\caO}{\mathcal{O}}
\newcommand{\caP}{\mathcal{P}}

\newcommand{\caR}{\mathcal{R}}
\newcommand{\caS}{\mathcal{S}}

\newcommand{\caW}{\mathcal{W}}

\newcommand{\field}{\mathbb}

\newcommand{\ga}{\alpha}
\newcommand{\gb}{\beta}

\newcommand{\C}{{\field C}}

\newcommand{\N}{{\field N}}
\newcommand{\Q}{{\field Q}}

\newcommand{\zz}{{\mathbf{z}}}

\theoremstyle{plain}

\newtheorem*{Fact}{Fact}

\newcommand{\excise}[1]{}

\newcommand{\cellsize}{11}
\newlength{\cellsz} 
\newsavebox{\cell}
\sbox{\cell}{\begin{picture}(\cellsize,\cellsize)
\put(0,0){\line(1,0){\cellsize}}
\put(0,0){\line(0,1){\cellsize}}
\put(\cellsize,0){\line(0,1){\cellsize}}
\put(0,\cellsize){\line(1,0){\cellsize}}
\end{picture}}
\newcommand\cellify[1]{\def\thearg{#1}\def\nothing{}%
\ifx\thearg\nothing
\vrule width0pt height\cellsz depth0pt\else
\hbox to 0pt{\usebox{\cell} \hss}\fi%
\vbox to \cellsz{
\vss
\hbox to \cellsz{\hss$#1$\hss}
\vss}}
\newcommand\tableau[1]{\vtop{\let\\\cr
\baselineskip -16000pt \lineskiplimit 16000pt \lineskip 0pt
\ialign{&\cellify{##}\cr#1\crcr}}}
\newcommand{\kellsize}{24}
\newlength{\kellsz} 
\newsavebox{\kell}
\sbox{\kell}{\begin{picture}(\kellsize,\kellsize)
\put(0,0){\line(1,0){\kellsize}}
\put(0,0){\line(0,1){\kellsize}}
\put(\kellsize,0){\line(0,1){\kellsize}}
\put(0,\kellsize){\line(1,0){\kellsize}}
\end{picture}}
\newcommand\kellify[1]{\def\thearg{#1}\def\nothing{}%
\ifx\thearg\nothing
\vrule width0pt height\kellsz depth0pt\else
\hbox to 0pt{\usebox{\kell} \hss}\fi%
\vbox to \kellsz{
\vss
\hbox to \kellsz{\hss$#1$\hss}
\vss}}
\newcommand\ktableau[1]{\vtop{\let\\\cr
\baselineskip -16000pt \lineskiplimit 16000pt \lineskip 0pt
\ialign{&\kellify{##}\cr#1\crcr}}}

\newcommand{\sellsize}{63}
\newlength{\sellsz} 
\newsavebox{\sell}
\sbox{\sell}{\begin{picture}(\sellsize,20)
\put(0,0){\line(1,0){\sellsize}}
\put(0,0){\line(0,1){\sellsize}}
\put(\sellsize,0){\line(0,1){\sellsize}}
\put(0,\sellsize){\line(1,0){\sellsize}}
\end{picture}}
\newcommand\sellify[1]{\def\thearg{#1}\def\nothing{}%
\ifx\thearg\nothing
\vrule width0pt height\sellsz depth0pt\else
\hbox to 0pt{\usebox{\sell} \hss}\fi%
\vbox to \sellsz{
\vss
\hbox to \sellsz{\hss$#1$\hss}
\vss}}
\newcommand\stableau[1]{\vtop{\let\\\cr
\baselineskip -16000pt \lineskiplimit 16000pt \lineskip 0pt
\ialign{&\sellify{##}\cr#1\crcr}}}

\begin{document}
\pagestyle{plain}
\mbox{}
\title{Governing singularities of symmetric orbit closures}
\author{Alexander Woo}
\author{Benjamin J.~Wyser}
\author{Alexander Yong}
\address{Dept~of Mathematics, University of Idaho, Moscow, ID 83844}
\email{awoo@uidaho.edu}
\address{Dept~of Mathematics, University of Illinois at
Urbana-Champaign, Urbana, IL 61801, USA}
\email{bwyser@uiuc.edu}
\address{Dept~of Mathematics, University of Illinois at
Urbana-Champaign, Urbana, IL 61801, USA}
\email{ayong@uiuc.edu}
\date{\today}
\maketitle
\begin{abstract}
We develop \emph{interval pattern avoidance} and \emph{Mars-Springer
ideals} to study  singularities of symmetric orbit closures in a flag variety. This paper focuses on the case 
of the Levi subgroup $GL_p\times GL_q$ acting on the classical flag variety.
We prove that all reasonable singularity properties can be classified in terms
of interval patterns of \emph{clans}.

\end{abstract}

\section{Introduction}
\subsection{Overview}

Let $G/B$ be a generalized flag variety where $G$
is a complex, reductive algebraic group
and $B$ is a choice of Borel subgroup.
A subgroup $K$
of $G$ is {\bf symmetric} if $K=G^{\theta}$
is the fixed point subgroup for an involutive automorphism 
$\theta$ of $G$.  Such a subgroup is {\bf spherical}, which means
that its action on $G/B$ by left translations has finitely many orbits.

The study of orbits of a symmetric subgroup on the flag variety
was initiated in work of G.~Lusztig--D.~Vogan~\cite{LV-83, Vogan-83}, who related the singularities
of their closures to characters of particular infinite-dimensional
representations of a certain real form $G_{\mathbb R}$ of $G$.
Since then, there has been a stream of results on the combinatorics
and geometry of these orbit closures.  Notably,
R.~W.~Richardson--T.~A.~Springer~\cite{Richardson-Springer-90}
gave a description of the partial order given by inclusions of orbit
closures, and M.~Brion~\cite{Brion-01} studied general properties of their
singularities, showing that, in many cases, including the one addressed
in this paper, all these orbit closures are normal and Cohen-Macaulay with
rational singularities.  One might also hope that the study of
singularities on closures of symmetric subgroup orbits would lead to
better understanding of the general relationship between the combinatorics
associated to spherical varieties and singularities of orbit closures on them; 
N.~Perrin has written a survey
of this topic \cite{Perrin-12}.

We initiate a combinatorial approach, backed by explicit commutative algebra computations, to the study of the singularities of these orbit closures. This paper
considers the case of the symmetric subgroup $K=GL_p \times GL_q$ of block diagonal matrices in $GL_n$, where $n=p+q$. Here $\theta$ is defined by
\begin{equation}
\label{eqn:theinvolnov29}
\theta(M)=I_{p,q}M I_{p,q},
\end{equation}
where
$I_{p,q}$ is the diagonal $\pm 1$ matrix with $p$ many $1$'s followed by $q$ many $-1$'s on the diagonal.  In this case, $G_{\mathbb R}=U(p,q)$, the \emph{indefinite unitary
group} of signature $(p,q)$.

For this symmetric subgroup $K$, the finitely many $K$-orbits ${\mathcal O}_{\gamma}$ and 
their closures $Y_{\gamma}$ can be parameterized by
{\bf $(p,q)$-clans} $\gamma$ \cite[Theorem 4.1]{Matsuki-Oshima-90} (see also \cite[Theorem 2.2.8]{Yamamoto-97}).  These
clans are partial matchings of vertices $\{1,2,\ldots,n\}$ where
each unmatched vertex is assigned a sign of $+$ or $-$; the difference in the number of $+$'s and $-$'s must be $p-q$. We represent clans by
length $n$ strings in ${\mathbb N}\cup \{+,-\}$, with pairs of equal numbers indicating a matching. Let ${\tt Clans}_{p,q}$ denote the set of all such clans.  For example, three clans from ${\tt Clans}_{7,3}$ are represented by the strings:
\[1+++2+1-2+, \ \  \ +1-++1+2+2, \ \text{\ \ and \ \ } +\!+\!-\!+\!+\!+\!-\!+\!+\!-.\]
Note that strings that differ by a permutation of the natural numbers, such as
$1212$ and $2121$, represent the same clan.

One inspiration for this work is W.~M.~McGovern's
characterization of singular $K$-orbit closures in terms of \emph{pattern avoidance} of clans \cite{McGovern-09a}.  
Suppose $\gamma \in {\tt Clans}_{p,q}$ and $\theta\in {\tt Clans}_{r,s}$.  Then $\theta=\theta_1 \ldots \theta_{r+s}$
is said to {\bf (pattern) include} $\gamma=\gamma_1 \ldots \gamma_{p+q}$ if there are
indices $i_1<i_2<\cdots<i_{p+q}$ such that:
\begin{enumerate}
 \item if $\gamma_j = \pm$ then $\theta_{i_j}=\gamma_j$; and
 \item if $\gamma_k=\gamma_{\ell}$ then $\theta_{i_k}=\theta_{i_{\ell}}$.
\end{enumerate}
For example, the clan $\gamma=1++--1$ contains the pattern $\theta=1++-1$, taking $(i_1,\hdots,i_5)$ to be either
$(1,2,3,4,6)$ or $(1,2,3,5,6)$. 

Say that $\theta$ {\bf (pattern) avoids} $\gamma$ if $\theta$ does not include $\gamma$.
The main theorem of \cite{McGovern-09a} asserts that $Y_{\gamma}$ is smooth if and only if $\gamma$ avoids
the patterns
\[1+-1, \ 1-+1, \ 1212, \ 1+221, \ 1-221, \ 122+1, \ 122-1, \ 122331.\]

On the other hand, in \cite[Section~3.3]{Woo-Wyser-14} it is noted
that $Y_{1++-1}$ is non-Gorenstein, while $Y_{1++--1}$ is Gorenstein, even though $1++--1$ pattern
contains $1++-1$. Therefore, a more general notion will sometimes be required to characterize
which $K$-orbits satisfy a particular singularity property.

Suppose that ${\mathcal P}$ is any {\bf singularity mildness property},
by which we mean a local property of
varieties that holds on open subsets and is stable under smooth morphisms.
Many singularity properties, such as being Gorenstein,
being a local complete intersection (lci), being factorial, having Cohen-Macaulay rank 
$\leq k$, or having Hilbert-Samuel multiplicity $\leq k$,
 satisfy these conditions. For such a $\mathcal{P}$, consider two related problems:
\begin{itemize}
\item[(I)] Which $K$-orbit closures $Y_{\gamma}$ are globally ${\mathcal P}$?
\item[(II)] What is the non-${\mathcal P}$-locus of $Y_{\gamma}$?
\end{itemize}

This paper gives a universal combinatorial language, {\bf interval pattern avoidance} of clans,
to answer these questions for any singularity mildness property, at the cost of
potentially requiring an infinite number of patterns.  This language is also useful for
collecting and analyzing data and partial results. We present explicit equations
for computing whether a property holds at a specific orbit $\caO_\alpha$ on an orbit
closure $Y_\gamma$.  Note that all points of ${\mathcal O}_{\alpha} \subseteq Y_{\gamma}$ are locally
isomorphic to one another, since the $K$-action can be used to move
any point of $\caO_{\alpha}$ isomorphically to any other point.
Since the non-$\caP$-locus is closed, it is a union of $K$-orbit closures.
Consequently, for any given clan $\gamma$, (II) can be answered
by finding a finite set of clans $\{\alpha\}$, namely those
indexing the irreducible components of the non-$\caP$-locus.
(I) asks if this set is nonempty.

This situation parallels that for Schubert varieties. 
In that setting, the first and third authors introduced interval pattern avoidance for permutations, showing that it provides a common
perspective to study all reasonable singularity measures~\cite{Woo-Yong:Gorenstein, Woo-Yong-08, Woo-2010}. 
 This paper gives the first analogue of those results for $K$-orbit closures.

Some properties ${\mathcal P}$ hold globally on every $Y_{\gamma}$. For those cases, the above questions are unnecessary.  For example, this is true when ${\mathcal P}=$``normal'' and 
${\mathcal P}=$``Cohen-Macaulay'' in the case $(G,K)=(GL_n,GL_p \times GL_q)$ of this paper.
(This is not the case for all symmetric pairs $(G,K)$.)

As is explained in \cite{McGovern-09a}, the above
answer to (I) for ${\mathcal P}=$``smooth'' 
also is the answer for the property $\caP=$``rationally smooth''.
(Recall rational smoothness 
means that the local intersection cohomology of $Y_\gamma$ (at a point of
$\caO_\alpha$) is trivial.)  However, the answer to question (II), which asks for a combinatorial description of the (rationally) singular locus, is unsolved except in some special cases~\cite{Woo-Wyser-14}. Actually, it is unknown whether the singular locus and rationally singular locus coincide
for all orbit closures (but see Conjecture \ref{conj:ratl-smooth-locus}).

For most finer singularity mildness properties, answers to both (I) and (II) are unknown. Perhaps the most famous such property
comes from the \emph{Kazhdan-Lusztig-Vogan} (KLV) polynomials $P_{\gamma,\alpha}(q)$. These polynomials are the link to representation theory that originally motivated the study of $K$-orbit closures.  For $(G,K)=(GL_n, GL_p\times GL_q)$, the KLV polynomial $P_{\gamma,\alpha}(q)\in {\mathbb Z}_{\geq 0}[q]$
is the Poincar\'e polynomial for the local intersection cohomology of $Y_{\gamma}$ at any point of the orbit ${\mathcal O}_{\alpha}$ \cite{LV-83}. Rational smoothness of $Y_{\gamma}$ along ${\mathcal O}_\alpha$ is hence equivalent to 
the equality $P_{\gamma,\alpha}(q)=1$.  More generally, for any fixed $k>1$, the property ${\mathcal P}=$``$P_{\gamma,\alpha}(1)\leq k$''
behaves as a singularity mildness property on $K$-orbit closures by recent work of W.~M.~McGovern \cite{McGovern:uppersemi}, but his proof uses representation theory, and the algebraic geometry is not well-understood.

\subsection{Main ideas}
\label{sec:main-results}
Each $K$-orbit closure $Y_{\gamma}$ is a union of $K$-orbits ${\mathcal O}_{\alpha}$; 
let $\leq$ denote the {\bf Bruhat (closure) order} on $K$-orbit closures on $GL_n/B$.  (This means $\alpha \leq \gamma$ if and
only if $Y_{\alpha} \subseteq Y_{\gamma}$). Also, let $[\alpha,\gamma]$ and $[\beta,\theta]$ be intervals in Bruhat order on ${\tt Clans}_{p,q}$
and ${\tt Clans}_{r,s}$, respectively. 

Define $[\beta,\theta]$ to {\bf interval pattern contain}
$[\alpha,\gamma]$, and write $[\alpha,\gamma]\hookrightarrow
[\beta,\theta]$ if: 
\begin{itemize}
\item[(a)] there are indices 
\[I:i_1<i_2<\cdots<i_{p+q}\] which commonly witness the containment of
$\gamma$ into $\theta$ and $\alpha$ into $\beta$, 
\item[(b)] $\theta$ and $\beta$ agree outside of these indices; and
\item[(c)] $
\ell(\theta)-\ell(\beta)=\ell(\gamma)-\ell(\alpha)$,
\end{itemize}
where if $\gamma=\gamma_1 \gamma_2 \ldots \gamma_{p+q}$ then
\[\ell(\gamma):=\sum_{\gamma_i=\gamma_j \in \N} j-i-\#\{\gamma_s=\gamma_t: s<i<t<j\}.\]
 
Notice $\beta$ is determined by $\alpha,\gamma$, $\theta$, and the set of indices $I$.
In particular, $\beta$ is the unique clan $\Phi(\alpha)$ that agrees with $\alpha$ on $I$ and agrees with $\theta$ on
$\{1,\ldots,n\}\setminus I$. Thus, we define a clan $\theta$ to {\bf interval pattern contain} $[\alpha,\gamma]$ if
$[\alpha,\gamma]$ is contained in $[\Phi(\ga),\theta]$. Similarly, we can speak of $\theta$ avoiding a list
of intervals.

\begin{Example}
\label{exa:intro1}
Let $[\alpha,\gamma]=[+--+, 1212]$ and $[\beta,\theta]=[1+--+1, 123231]$. Then one can check $[\beta,\theta]$ interval contains $[\alpha,\gamma]$, using the middle four positions. \qed
\end{Example}
\begin{Example}
\label{exa:intro2}
Let $[\alpha,\gamma]=[+-, 11]$ and $[\beta,\theta]=[\underline{+}+\underline{-}, \underline{1}+\underline{1}]$.
If $[\beta,\theta]$ contains $[\alpha,\gamma]$ it must be using the underlined positions. However, $\ell(\gamma)-\ell(\alpha)=1$ while $\ell(\theta)-\ell(\beta)=2$.  Thus $[\gb,\theta]$ does not contain $[\ga,\gamma]$.\qed
\end{Example}

Define
\[{\mathcal C}:=\{[\alpha,\gamma]: \alpha \leq \gamma \mbox{\ in some \ } {\tt Clans}_{p,q}\}\subseteq {\tt Clans}\times {\tt Clans}\]
where 
\[{\tt Clans}=\union_{p,q}{\tt Clans}_{p,q}.\] 
Declare $\preceq_{\mathcal C}$ to be the poset relation on ${\mathcal C}$ generated by
\begin{itemize}
\item $[\alpha,\gamma]\preceq_{{\mathcal C}}[\beta,\theta]$ if $[\beta,\theta]$ interval pattern contains
$[\alpha,\gamma]$, and
\item $[\alpha,\gamma]\preceq_{{\mathcal C}}[\alpha',\gamma]$ if $\alpha' \leq \alpha$.
\end{itemize}

For any poset $(\caS,\preceq)$, an \textbf{upper order ideal} is a subset $\caI$ of $\caS$ having the property
that whenever $x \in \caI$, we also have $y \in \caI$ for all $y \geq x$. The following theorem
provides a basic language to express 
answers to (I) and (II). 

\begin{Theorem}
\label{thm:main}
Let ${\mathcal P}$ be a singularity mildness property (see Definition~\ref{def:mildness}).
\begin{itemize}
\item[(I)] The set of intervals $[\alpha,\gamma] \in {\mathcal C}$
such that ${\mathcal P}$ fails on each point of ${\mathcal O}_{\alpha}\subseteq Y_{\gamma}$ is an upper order ideal in $(\mathcal{C}, \preceq_{\mathcal{C}})$.
\item[(II)]  The set of clans $\gamma$ such
that ${\mathcal P}$ holds at all points of $Y_{\gamma}$ are those that avoid all the intervals
$[\alpha_i, \gamma_i]$ constituting some (possibly infinite) set ${\mathcal A}_{\mathcal P} \subseteq {\mathcal C}$.
\end{itemize}
\end{Theorem}

For simplicity, we work over $\mathbb{C}$, but our results are valid over any field of characteristic not 2, with the caveat that
even for the same $\caP$, the set $\caA_{\caP}$ may be
field dependent.

Although stated in combinatorial language, as we will see, Theorem~\ref{thm:main} follows from a geometric result, 
Theorem~\ref{thm:iso}, which
establishes a local isomorphism between certain ``slices'' of the orbit closures.  Rather than working directly with the
$K$-orbit closures on $G/B$, it is easier to establish
this isomorphism using particular slices of $B$-orbit closures on $G/K$ given by J.~G.~M.~Mars and
T.~A.~Springer \cite{Mars-Springer}.  

The space $G/K=GL_{n}/(GL_p\times GL_q)$ is
the configuration space of all splittings of $\C^n$ as a direct sum
$V_1 \oplus V_2$ of subspaces, with $\dim(V_1) = p$ and $\dim(V_2)=q$.  Indeed, $G=GL_n(\C)$ acts transitively on
the set of all such splittings, and $K=GL_p\times 
GL_q$ is the stabilizer of the ``standard'' splitting of
$\C^n$ as $\left\langle e_1,\hdots,e_p \right\rangle \oplus  \left\langle e_{p+1},\hdots,e_n \right\rangle$, where
$e_1,\hdots,e_n$ is the standard basis of $\C^n$.  A point $gK/K \in G/K$, with $g \in GL_n(\C)$ an invertible $n \times n$
matrix, is identified with the splitting whose $p$-dimensional component is the span of the first $p$ columns 
of $g$ and whose $q$-dimensional component is the span of the last $q$ columns of $g$.

$B$-orbits on $G/K$ are in bijection with $K$-orbits on
$G/B$ (both orbit sets being in bijection with the $B \times K$-orbits on $G$). In addition, this bijection
preserves all mildness properties.  This is because the orbit
closures correspond via the two locally trivial fibrations
$G \rightarrow G/K$ and $G \rightarrow G/B$, each of which has smooth fiber.
The $B$-orbit on $G/K$ corresponding
to the clan $\gamma$ will be
denoted by $Q_{\gamma}$, and its closure will be denoted by $W_{\gamma}$.  For any singularity mildness property $\mathcal{P}$, this discussion implies the following:

\begin{Observation}\label{obs:orbits-equivalent}
 $Y_{\gamma}$ is $\mathcal{P}$ along $\caO_{\ga}$ if and only if $W_{\gamma}$ is $\mathcal{P}$ along $Q_{\ga}$.
\end{Observation}

In \cite{Wyser-Yong-13}, the second and third authors considered certain open affine subsets of $K$-orbit closures
on $G/B$ which they called \emph{patches}.  That work also introduced \emph{patch ideals} of equations which set-theoretically cut out the
patches.  Experimental computation using patch ideals led to several of the conjectures which appear in
\cite{Wyser-Yong-13, Woo-Wyser-14}.

In contrast, the Mars-Springer slices are not open affine pieces of the orbit closures.
However, as we will see, they are essentially as good, in that they carry all of the local information
that we are interested in.  They are analogues of \emph{Kazhdan-Lusztig varieties}, which play a similar role in the
study of Schubert varieties.

In this paper we introduce \emph{Mars-Springer ideals} (see Section~\ref{sec:ms-variety})
with explicit equations that
set-theoretically cut out the Mars-Springer slices.
These equations are conjectured to also be scheme-theoretically correct; see Conjecture~\ref{conj:radical}.
These are the $K$-orbit versions of \emph{Kazhdan-Lusztig ideals}, which define the aforementioned Kazhdan-Lusztig varieties.
The latter ideals have
been of use in both computational and theoretical analysis of Schubert varieties (see, for example,
\cite{Woo-Yong-Grobner-12, Ulfarsson-Woo-11} and the references
therein).  In the same vein, we mention a practical advantage of the Mars-Springer ideal over the patch ideal.  Gr\"{o}bner basis calculations with the former are
several times faster than those of the latter, as fewer variables are involved.

\subsection{Organization} In Section~2, we present some preliminaries. In particular, we more precisely define
``mildness properties,'' giving examples and establishing some basic facts that we will need.  We next recall the
attractive slices of Mars-Springer \cite{Mars-Springer}. Finally, we describe the
quasi-projective variety structure of $G/K$ for the case $(G,K)=(GL_n,GL_p \times GL_q)$ of this paper. 
In Section~3, we give explicit affine coordinates for the Mars-Springer slice (Theorem~\ref{thm:MS-variety-is-slice}).
Section~4 defines the \textit{Mars-Springer variety} and its ideal. It culminates with Theorem~\ref{thm:iso}, which asserts that certain Mars-Springer varieties are
isomorphic to others. This theorem is the key to our proof of Theorem~\ref{thm:main}. 
In order to prove Theorem~\ref{thm:iso},
we develop combinatorics of interval embeddings using earlier work of the second author
\cite{Wyser-15} on the Bruhat order on clans; this is Section 5.
We then give our proofs of Theorem~\ref{thm:iso} and Theorem~\ref{thm:main}
in Section~6.
We conclude with problems and conjectures in Section~7. 

\section{Preliminaries}
\subsection{Singularity Mildness Properties} 
We define the class of local properties that we are interested in, with specific examples.

\begin{Definition}\label{def:mildness}
Suppose that $\caP$ is a local property of algebraic varieties, meaning that $\caP$ is verified at a point $w$ of an
algebraic variety $W$ based solely on the local ring $\caO_{W,w}$.  We say that $\caP$ is a \textbf{singularity mildness property} (or simply a \textbf{mildness property}) if $\caP$ is
\begin{enumerate}
 \item Open, meaning that the $\caP$-locus of any variety is open; and
 \item Stable under smooth morphisms.  Precisely, we mean that if $X$ and $Y$ are any varieties, and $f: X \rightarrow Y$
 is any smooth morphism between them, then for any $x \in X$, $X$ is $\caP$ at $x$ if and only if $Y$ is $\caP$ at $f(x)$.
\end{enumerate}
\end{Definition}

For any smooth variety $S$ and for any variety $X$, the projection map $S \times X \rightarrow X$ is a smooth morphism.
Thus:

\begin{Observation}\label{obs:proj-smooth}
 If $\caP$ is a mildness property, $S$ is a smooth variety, and $X$ is any variety, then for any $x \in X$ and 
 for any $s \in S$, $X$ is $\caP$ at $x$ if and only if $S \times X$ is $\caP$ at $(s,x)$.
\end{Observation}

The next lemma gives a (non-exhaustive) list of examples of mildness properties $\mathcal{P}$.

\begin{Lemma}\label{lem:slice-stable-properties}
 Examples of mildness properties $\mathcal{P}$ include:
 \begin{enumerate}[(i)]
  \item Reducedness;
  \item Normality;
  \item Smoothness;
  \item lci-ness;
  \item Gorensteinness.
 \end{enumerate}
\end{Lemma}
\begin{proof}
That these properties are open (on varieties) is well-known.  That they are stable under smooth morphisms
follows from results of \cite[\S 23]{Matsumura}.  Indeed, let $X$ and $Y$ be varieties, with $f: X \rightarrow Y$
a smooth morphism between them.  Let $x \in X$ be given, and let $y=f(x)$.  Now, $f$ is flat, and also the fiber
$F_y$ over $y$ is smooth over the residue field $\kappa(y)$, hence regular.  Thus $\mathcal{O}_{F_y,x}$ is regular,
hence also lci and Gorenstein.  Now (iii) and (v) follow from Theorems 23.7 and 23.4 \cite[\S 23]{Matsumura}, respectively,
while (iv) is mentioned in the remark following Theorem 23.6 of \emph{loc. cit.}; as indicated there, further details
can be found in \cite{Avramov}.

For (i) and (ii), we appeal to \cite[Corollary to Theorem 23.9]{Matsumura}.  Note that this result requires that the fiber
ring corresponding to any prime ideal of $\mathcal{O}_{Y,y}$ be reduced (resp. normal), rather than just the fiber ring
of the maximal ideal.  In our case, all such fiber rings are once again regular (hence both reduced and normal).  Indeed,
if $p$ is any prime of $\mathcal{O}_{Y,y}$, then $p$ corresponds to an irreducible closed subvariety $Y'$ of $Y$.  The
inverse image scheme $f^{-1}(Y')$ is smooth over $Y'$, and the fiber ring
$\mathcal{O}_{X,x} \otimes_{\mathcal{O}_{Y,y}} \kappa(p)$ is the generic fiber of $f^{-1}(Y') \rightarrow Y'$.
\end{proof}

\subsection{Mars-Springer slices}\label{sec:attractive-slices}
Let $W$ be any irreducible variety with an action of a Borel group $B$, let $w \in W$ be given, and
let $Q=B \cdot w$.  The following definition is from \cite{Mars-Springer}:

\begin{Definition}\label{def:slice}
An \textbf{attractive slice} to $Q$ at $w$ in $W$
is a locally closed subset $S$ of $W$ containing $w$ such that:
\begin{enumerate}
	\item The restriction of the action map $B \times S \rightarrow W$ is a smooth morphism;
	\item $\dim(S)+\dim(Q)=\dim(W)$; and
	\item There exists a map $\lambda: \mathbb{G}_m \rightarrow B$ such that
	  \begin{enumerate}
	    \item $S$ is stable under $Im(\lambda)$;
	    \item $w$ is a fixed point of $Im(\lambda)$; and
	    \item The action map $\mathbb{G}_m \times S \rightarrow S$ extends to a morphism
		  $\mathbb{A}^1 \times S \rightarrow S$ which sends $\{0\} \times S$ to $w$.
	  \end{enumerate}
\end{enumerate}
\end{Definition}

The following result justifies our use of slices in studying singularities, since taking slices preserves 
the local properties that we are interested in.
\begin{Lemma}\label{lem:mildness-property-is-slice-stable}
 Given notation as in Definition~\ref{def:slice} and a mildness property $\caP$, the variety $W$ is $\caP$ at $w$
 if and only if $S$ is $\caP$ at $w$.
\end{Lemma}
\begin{proof}
 Indeed, the action map $B \times S \rightarrow W$ is smooth, so $W$ is $\caP$ at $w$ if and only if $B \times S$ is
 $\caP$ at $(1,w)$, which is the case if and only if $S$ is $\caP$ at $w$, by Observation~\ref{obs:proj-smooth}.
\end{proof}

Let $G$ be a complex algebraic group with an involutive automorphism $\theta$, with $K=G^{\theta}$ the corresponding symmetric
subgroup.  Suppose that $T \subseteq B$ are a $\theta$-stable maximal torus and Borel subgroup
of $G$, respectively.  Let $N$ be the normalizer of $T$ in $G$, and let $\caW = N/T$ be the Weyl group.

Given a $B$-orbit $Q$ on $G/K$, we now define the {\bf Mars-Springer slice} $S_{\overline{x}}$ to 
$Q$ at a specially chosen point $\overline{x}\in Q$.  Letting $U^-$ be the unipotent radical of the Borel subgroup
opposite to $B$, this slice is of the form
\begin{equation}\label{eqn:slice}
 S_{\overline{x}} := (U^- \cap \psi(U^-)) \cdot \overline{x},
\end{equation}
where $\psi$ is a certain involution on $G$.

To be precise, consider the set
\[\mathcal{V} := \{x \in G \mid x(\theta x)^{-1} \in N\}.\]  

In \eqref{eqn:slice}, we pick any 
\[\overline{x}:=xK/K\in Q \text{\ with $x\in {\mathcal V}$;}\]
such a choice exists by \cite[Theorem 1.3]{Richardson-Springer-90}.
Let $\eta$ be the image in $\caW$ of $x \theta(x)^{-1} \in N$.  The element $x \theta(x)^{-1} \in N$ may depend upon the
choice of $\overline{x}$, but $\eta$, its class modulo $T$, does not.  Then in (\ref{eqn:slice}) we take 
\[ \psi := c_{\eta} \circ \theta, \]
where $c_\eta$ denotes conjugation by $\eta$. With these choices, 
it is shown in \cite[Section 6.4]{Mars-Springer} that
$S_{\overline{x}}$, as defined in \eqref{eqn:slice}, is indeed an attractive slice to $Q$ at $\overline{x}$ in
$G/K$, in the sense of Definition~\ref{def:slice}.

The following lemma is more or less formal; we include a proof for completeness:

\begin{Lemma}\label{lem:slice-of-orbit-closure}
With notation as above, if $W$ is any $B$-orbit closure on $G/K$ containing $Q$, then 
\begin{enumerate}[(i)]
 \item The scheme-theoretic intersection $W \cap S_{\overline{x}}$ is reduced, and
 \item $W \cap S_{\overline{x}}$ is an attractive slice to $Q$ at $\overline{x}$ in $W$.
\end{enumerate}
\end{Lemma}
\begin{proof}
It is easy to see that the diagram
\[
\xymatrix{
  B \times (W \cap S_{\overline{x}}) \ar[r] \ar@{^{(}->}[d] & W \ar@{^{(}->}[d] \\
  B \times S_{\overline{x}} \ar[r] & G/K
}
\]
is Cartesian.  Therefore, since the action map along the bottom is smooth, the (restricted) action map
along the top is also smooth, since smooth morphisms are stable under base extension.  Since $W$ is reduced and
reducedness is a mildness property by Lemma~\ref{lem:slice-stable-properties}, $B \times (W \cap S_{\overline{x}})$
is reduced, and thus $W \cap S_{\overline{x}}$ is reduced, by Observation~\ref{obs:proj-smooth}.  This proves (i).

For (ii), the preceding observation regarding the smoothness of the map along the top of the above diagram
verifies Definition~\ref{def:slice}(1) for $S=W \cap S_{\overline x}$.  This map furthermore
has the same relative dimension as the one along the bottom, which implies that
\[ \dim(G/K) - \dim(W) = \dim(S_{\overline{x}}) - \dim(W \cap S_{\overline{x}}). \]
Since $S_{\overline{x}}$ is an attractive slice to $Q$ in $G/K$, we know that 
\[\dim(S_{\overline{x}}) + \dim(Q) = \dim(G/K),\] 
and combining these two equations gives 
\[\dim(Q) + \dim(W \cap S_{\overline{x}}) = \dim(W),\] 
as required by Definition~\ref{def:slice}(2).

That $W \cap S_{\overline{x}}$ satisfies Definition \ref{def:slice}(3)
is obvious:  Let $\lambda: \mathbb{G}_m \rightarrow B$ be a one-parameter subgroup having properties (a) and (b)
relative to the slice $S_{\overline{x}}$.  Then clearly $W \cap S_{\overline{x}}$ is stable under $Im(\lambda)$,
$\overline{x}$ is still a fixed point of $Im(\lambda)$, and the action map
$\mathbb{G}_m \times (W \cap S_{\overline{x}}) \rightarrow W \cap S_{\overline{x}}$ still extends to a map
$\mathbb{A}^1 \times (W \cap S_{\overline{x}}) \rightarrow W \cap S_{\overline{x}}$ sending
$\{0\} \times (W \cap S_{\overline{x}})$ to $\overline{x}$ (simply restrict the original extension to the smaller subset).
\end{proof}


\begin{Corollary}\label{cor:irreducible}
 For $(G,K)=(GL_n,GL_p \times GL_q)$, let clans $\ga$ and $\gamma$ be given, with $\ga \leq \gamma$ in Bruhat order.
 Let $X_{\ga}$ denote the Mars-Springer slice to $Q_{\ga}$ at a particular point $x_{\ga} \in Q_{\ga}$ in $G/K$.  Then
 the slice $W_{\gamma} \cap X_{\ga}$ is irreducible.
\end{Corollary}
\begin{proof}
 It is observed in \cite{Brion-01} that all $B$-orbit closures on $G/K$
 for this particular case are ``multiplicity-free'', which implies by \cite[Theorem 5]{Brion-01} that $W_{\gamma}$ has
 rational singularities.  In particular, it is normal.  Since normality is a mildness property by Lemma
 \ref{lem:slice-stable-properties}, $W_{\gamma} \cap X_{\ga}$ is also normal.  Now, the ``attractive'' condition of
 slices (cf. part (3) of Definition \ref{def:slice}) ensures that $W_{\gamma} \cap X_{\ga}$ is connected because $x_{\ga}$
 must lie in every connected component of it.  Being both normal and connected, $W_{\gamma} \cap X_{\ga}$ must be
 irreducible.
\end{proof}

\begin{Remark}
 The result of Corollary \ref{cor:irreducible} in fact holds for any symmetric pair $(G,K)$ such that
 all of the $B$-orbit closures on $G/K$ are multiplicity-free.  Other such cases include $(GL_{2n},Sp_{2n})$ and
 $(SO_{2n},GL_n)$, as noted in \cite{Brion-01}.  However, not all symmetric pairs $(G,K)$ have this property.
 For those that do not, $B$-orbit closures on $G/K$ can be non-normal, and the
 Mars-Springer slices can in fact be reducible.  One example is $(GL_n,O_n)$.
 \qed
\end{Remark}

\subsection{$GL_{p+q}/(GL_p\times GL_q)$ as a quasi-projective 
variety}\label{sec:G/K-patches}

 We now discuss the structure of $G/K=GL_{p+q}/(GL_p\times GL_q)$ as a quasi-projective variety and its affine patches.
 The statements of this section
 could be extracted from standard definitions and results, but to the best of our knowledge, this has never been
 made explicit in the literature, so we take the opportunity to do so here.
 
 As stated in the Introduction, $G/K$ is the configuration space whose points correspond to splittings of $\mathbb{C}^n$ as
 a direct sum $V_1\oplus V_2$ of subspaces, with $\dim(V_1)=p$ and $\dim(V_2)=q$.  Hence it is natural to identify
 $G/K$ with a subset of the product of Grasmmanians $\mathrm{Gr}(p,{\mathbb C}^n)\times\mathrm{Gr}(q,{\mathbb C}^n)$.
 We can realize $\mathrm{Gr}(p,{\mathbb C}^n)\times\mathrm{Gr}(q,{\mathbb C}^n)$ as a subvariety of
 $\mathbb{P}^{\binom{n}{p}-1}\times\mathbb{P}^{\binom{n}{q}-1}$ using the Pl\"ucker embedding into each factor.
 The Segre embedding then realizes this product as a subvariety of $\mathbb{P}^{\binom{n}{p}\binom{n}{q}-1}$, so
 $\mathrm{Gr}(p,\mathbb{C}^n)\times\mathrm{Gr}(q,\mathbb{C}^n)$ is a projective variety.
 
 We let $p_S$ and $q_T$ denote the Pl\"ucker coordinates on $\mathrm{Gr}(p,{\mathbb C}^n)$ and $\mathrm{Gr}(q,{\mathbb C}^n)$ respectively,
 where $S$ and $T$ respectively range over all subsets of $\{1,\ldots,n\}$ of size $p$ (respectively $q$).
 
 Not every pair $V_1$ and $V_2$ of subspaces of $\mathbb{C}^n$ gives a splitting. In order for $V_1$ and $V_2$
 to form a direct sum, a basis for $V_1$ must be linearly independent of a basis of $V_2$.  Let $M$ be a matrix
 whose first $p$ columns are a basis for $V_1$ and whose last $q$ columns are a basis for $V_2$.  If we take the
 Laplace expansion of $\det M$ using the first $p$ columns and identify determinants of submatrices with
 Pl\"ucker coordinates, then we get
 \begin{equation}
 \label{eqn:G/K-bdry}
 \det M=\sum_{S\subset\{1,\ldots,n\}, \#S=p} (-1)^S p_Sq_{\overline{S}},
 \end{equation}
 where $\overline{S}:=\{1,\ldots,n\}\setminus S$ and $(-1)^S:=(-1)^{(\sum_{s\in S} s)-\binom{p+1}{2}}$.
 Hence we can identify $G/K$ with the \emph{open} subset of $\mathrm{Gr}(p,{\mathbb C}^n)\times\mathrm{Gr}(q,{\mathbb C}^n)$ where
 (\ref{eqn:G/K-bdry}) is nonzero.  This gives $G/K$ the structure of a quasi-projective variety.
 
 The product $\mathbb{P}^{\binom{n}{p}-1}\times\mathbb{P}^{\binom{n}{q}-1}$ of projective spaces has a covering
 by affine patches $\{U_{S,T}:=U_S\times U_T\}$, where $U_S$ is the patch where $p_S$ is nonzero and $U_T$ is
 the patch where $q_T$ is nonzero.
 Thus $\mathrm{Gr}(p,\mathbb{C}^n)\times\mathrm{Gr}(q,\mathbb{C}^n)$ is covered by these affine patches.
 As explained in~\cite[Section 9.1]{Fulton-YoungTableaux}, any subspace $V_1 \in \mathrm{Gr}(p,\mathbb{C}^n)$ having
 Pl\"ucker coordinate $p_S(V_1) \neq 0$ has a basis $\{e_s+\sum_{r\not\in S} a_{s,r}(V_1)e_r\}_{s\in S}$ where $e_i$ denotes
 the $i$-th standard basis vector, and the functions $a_{s,r}$ are local coordinates on the patch
 $U_S \cap \mathrm{Gr}(p,\mathbb{C}^n)$.  Over this patch, there is an algebraic section $\phi_S$ of the projection
 map $M_{n \times p}^0 \to \mathrm{Gr}(p,\mathbb{C}^n)$, where $M_{n \times p}^0$ denotes the $n \times p$ matrices of
 full rank $p$.  The section $\phi$ sends the subspace $V_1$ to the matrix 
 whose columns are precisely the aforementioned basis.
 
 Likewise, any subspace
 $V_2 \in \mathrm{Gr}(q,\mathbb{C}^n)$ having
 Pl\"ucker coordinate $q_T(V_2) \neq 0$ has a basis $\{e_t + \sum_{r \not\in T} b_{t,r}(V_2)e_r\}_{t \in T}$,
 and the $b_{t,r}$ are local coordinates on the patch $U_T \cap \mathrm{Gr}(q,\mathbb{C}^n)$.  Over
 this patch, there is a similar section $\phi_T$ which sends $V_2$ to the full rank
 $n \times q$ matrix whose columns are this basis.
 
 Taking the products of the sections $\phi_S$ and $\phi_T$, over
 $U_{S,T} \cap (\mathrm{Gr}(p,\mathbb{C}^n)\times\mathrm{Gr}(q,\mathbb{C}^n))$ we have an algebraic map $\phi_{S,T}$ from $G/K$
 into $M_{n \times p}^0 \times M_{n \times q}^0$, which naturally embeds into $M_{n \times n}$ by simple concatenation
 of matrices.  Over $G/K$, this map takes values in $G=GL_n$ and gives an algebraic section $\phi_{S,T}$ of the
 projection map $G \to G/K$. 

\section{Affine coordinates for the Mars-Springer slice}
\label{sec:Mars-Springer-coordinates}
We provide explicit coordinates for the Mars-Springer slice for the case that $(G,K)=(GL_n,GL_p \times GL_q)$.

\subsection{The affine space $S_{\ga}$}\label{sec:affine-space-s}
We define an affine space of matrices associated to each clan $\alpha$. First,
let $w_{\alpha} \in S_n$ (in one line notation) be defined as follows.  From left to right, $1$ through $p$ are assigned to the $+$'s and
left ends of matchings. Assign $\{p+1,\hdots,n\}$ as we read the clan from left to right.  When we encounter a $-$, we assign the smallest unused number, and when we encounter the {\em left} end of a matching, we immediately assign the smallest unused number from $\{p+1,\hdots,n\}$ to the corresponding {\em right} end of the matching.

\begin{Example}
\label{exa:122133}
If $\alpha=122133\in {\tt Clans}_{3,3}$ then $w_{\alpha}=125436$.
\qed
\end{Example}

We now construct, in stages, the generic matrix $M_{\ga}(\mathbf{z})$ of $S_{\alpha}$.
First, for $i=1,\hdots,n$:
\begin{itemize}
	\item[(O.1)] If $\ga_i=\pm$, then row $i$ has a $1$ in position $w_{\ga}(i)$.
	\item[(O.2)] If $(i<j)$ is a matching of $\ga$, then row $i$ has $1$'s in positions
	      $w_{\ga}(i)$ and $w_{\ga}(j)$.
	\item[(O.3)] If $(j<i)$ is a matching of $\ga$, then row $i$ has a $-1$ in position
	      $w_{\ga}(j)$ and a $1$ in position $w_{\ga}(i)$.
\end{itemize}

A {\bf pivot} is the northmost $1$ in each column.  (Note that $w_{\ga}$ is chosen precisely so that the pivots in the first $p$
columns, and in the last $q$ columns, occur northwest to southeast.)  In the first $p$ columns only, 
set to $0$ all entries
\begin{itemize}
\item[(Z.1)] in the same row as a pivot;
\item[(Z.2)] above a pivot (and in that pivot's column);
\item[(Z.3)] between the pivot $1$ and the corresponding $-1$ of a column (for all columns to which this applies); or
\item[(Z.4)] right of a $-1$.
\end{itemize}

In the rightmost $q$ columns, set to $0$ all entries
\begin{itemize}
\item[(Z.5)] in the same row as a pivot;
\item[(Z.6)] above a pivot (and in that pivot's column);
\item[(Z.7)] between the pivot $1$ and the corresponding $1$ below it in the same column (for all columns to which this applies); or
\item[(Z.8)] right of a $1$ which is the second $1$ in its column.
\end{itemize}

The remaining entries $z_{ij}$ of $M_{\ga}(\mathbf{z})$ are arbitrary, with an exception for each pair of matchings
of $\ga$ which are ``in the pattern $1212$''.  For each such pair of matchings
$(i<k)$, $(j<\ell)$ with $i < j < k < \ell$:
\begin{itemize}
\item[(Z.9)] Position $(\ell, w_{\alpha}(i))$ has entry $z_{\ell, w_{\alpha}(i)}$
(as usual) whereas its ``partner'' position $(\ell, w_{\ga}(k))$ has entry $-z_{\ell, w_{\alpha}(i)}$.
\end{itemize}

\begin{Example}
Let $\ga=1+12-2\in {\tt Clans}_{3,3}$.
Then $w_{\ga}=124365$.  The $1$'s
and $-1$'s are first placed into a $6 \times 6$ matrix, followed by placement of $0$'s as follows:
\[
 \begin{pmatrix}
    1 & \cdot & \cdot & 1 & \cdot & \cdot \\
    \cdot & 1 & \cdot & \cdot & \cdot & \cdot \\
    -1 & \cdot & \cdot & 1 & \cdot & \cdot \\
    \cdot & \cdot & 1 & \cdot & 1 & \cdot \\
    \cdot & \cdot & \cdot & \cdot & \cdot & 1 \\
    \cdot & \cdot & -1 & \cdot & 1 & \cdot
 \end{pmatrix}
\mapsto
 \begin{pmatrix}
    1 & 0 & 0 & 1 & 0 & 0 \\
    0 & 1 & 0 & 0 & 0 & 0 \\
    -1 & 0 & 0 & 1 & 0 & 0 \\
    0 & 0 & 1 & 0 & 1 & 0 \\
    \cdot & \cdot & 0 & 0 & 0 & 1 \\
    \cdot & \cdot & -1 & \cdot & 1 & 0
 \end{pmatrix}
.
\]

Since $\ga$ has no $1212$-patterns, the remaining unspecialized entries of the matrix are arbitrary.
Thus $S_{1+12-2}$ is an affine space of dimension $5$ (the codimension of $Q_{1+12-2}$ in $GL_6/K$), with its generic
entry being the matrix
\[
M_{1+12-2}(\mathbf{z}) = \begin{pmatrix}
    1 & 0 & 0 & 1 & 0 & 0 \\
    0 & 1 & 0 & 0 & 0 & 0 \\
    -1 & 0 & 0 & 1 & 0 & 0 \\
    0 & 0 & 1 & 0 & 1 & 0 \\
    z_{5,1} & z_{5,2} & 0 & 0 & 0 & 1 \\
    z_{6,1} & z_{6,2} & -1 & z_{6,4} & 1 & 0
 \end{pmatrix}
.
\]
\qed
\end{Example}

\begin{Example}
Now, consider $\ga=1+21-2$.  Here $w_{\ga}=123465$.  The placement of $1$'s and $-1$'s, and then $0$'s, is achieved by:
\[
 \begin{pmatrix}
    1 & \cdot & \cdot & 1 & \cdot & \cdot \\
    \cdot & 1 & \cdot & \cdot & \cdot & \cdot \\
    \cdot & \cdot & 1 & \cdot & 1 & \cdot \\
    -1 & \cdot & \cdot & 1 & \cdot & \cdot \\
    \cdot & \cdot & \cdot & \cdot & \cdot & 1 \\
    \cdot & \cdot & -1 & \cdot & 1 & \cdot
 \end{pmatrix}\mapsto
 \begin{pmatrix}
    1 & 0 & 0 & 1 & 0 & 0 \\
    0 & 1 & 0 & 0 & 0 & 0 \\
    0 & 0 & 1 & 0 & 1 & 0 \\
    -1 & 0 & 0 & 1 & 0 & 0 \\
    \cdot & \cdot & 0 & 0 & 0 & 1 \\
    \cdot & \cdot & -1 & \cdot & 1 & 0
 \end{pmatrix}
.
\]
The $1212$-pattern in positions $1 < 3 < 4 < 6$ dictates that the unspecialized positions $(6,w_{\ga}(1))=(6,1)$ and
$(6,w_{\ga}(4))=(6,4)$ are negatives of one another.  Thus $S_{1+21-2}$ is an affine space of dimension $4$ (the
codimension of $Q_{1+21-2}$ in $GL_6/K$), and
\[
M_{1+21-2}(\mathbf{z}) = \begin{pmatrix}
    1 & 0 & 0 & 1 & 0 & 0 \\
    0 & 1 & 0 & 0 & 0 & 0 \\
    0 & 0 & 1 & 0 & 1 & 0 \\
    -1 & 0 & 0 & 1 & 0 & 0 \\
    z_{5,1} & z_{5,2} & 0 & 0 & 0 & 1 \\
    z_{6,1} & z_{6,2} & -1 & -z_{6,1} & 1 & 0
 \end{pmatrix}
.
\]
\qed
\end{Example}

\begin{Example}
Each $1212$-pattern gives rise to a separate
identification of coordinates. For instance, when $\ga=123123$, the reader can verify that 
\[
M_{123123}(\mathbf{z}) = \begin{pmatrix}
    1 & 0 & 0 & 1 & 0 & 0 \\
    0 & 1 & 0 & 0 & 1 & 0 \\
    0 & 0 & 1 & 0 & 0 & 1 \\
    -1 & 0 & 0 & 1 & 0 & 0 \\
    z_{5,1} & -1 & 0 & -z_{5,1} & 1 & 0 \\
    z_{6,1} & z_{6,2} & -1 & -z_{6,1} & -z_{6,2} & 1
 \end{pmatrix}
.
\]
\qed
\end{Example}

The following result will be used in Section~\ref{sec:ms-variety}.
\begin{Lemma}
\label{lemma:isinvertible}
$M_{\ga}(\zz)$ invertible over $\Q[\zz]$.
\end{Lemma}
\begin{proof}
By construction, the matrix is invertible over the field of rational functions ${\mathbb Q}(\zz)$. The lemma would follow if we could prove that the determinant is a constant (which will be necessarily nonzero).
We argue by induction on $n=p'+q'$ that $\det M_{\ga'}(\zz)\in {\mathbb Q}$ for $\ga'\in {\tt Clans}_{p',q'}$. The base case $n=1,2$ 
is easy to check.

    Suppose $\gamma\in {\tt Clans}_{p,q}$ where
    $n=p+q>2$. If $\gamma_1=+$ consider
    the submatrix of $M_{\ga}(\zz)$ obtained by striking the first
    row and column. One checks that this submatrix $M'$ is equal to $M_{\ga'}(\zz)$, where
    $\ga'$ is $\ga$ with the first $+$ removed. By induction, 
    $\det M'\in {\mathbb Q}$. Then by cofactor expansion along the first row, since the unique nonzero entry is a $1$ in position $(1,1)$, 
    we have that $\det M=\det M'\in {\mathbb Q}$, as desired. 
    
    A similar argument may be used if $\ga_1=-$.

Finally, suppose $\ga_1$ is the left end of a matching. Let $\ga_r$
be the right end of said matching. Observe that the submatrix
$M'$ obtained by deleting the first row and first column of $M=M_{\ga}$
is the matrix $M_{\ga'}$ where $\ga'$ is obtained by deleting $\ga_1$
and replacing $\ga_r$ by $-$. Thus by induction $\det M'\in {\mathbb Q}$. 

Similarly, consider the submatrix $M''$ obtained by deleting the first row
and $w(r)$-th column. We observe that by 
\begin{itemize}
\item cyclically reordering rows $1,2,\ldots,w(r)-1$ of $M''$; and then 
\item multiplying the first row by $-1$,
\end{itemize}
we obtain the matrix $M_{\ga''}$ where $\ga''$ is obtained
from $\ga$ by deleting $\ga_r$ and replacing $\ga_1$ by $+$.
Hence, by cofactor expansion along the first row of $M$, we obtain that
$\det M=\det M'\pm\det M''\in {\mathbb Q}$,
as desired.
\end{proof}

\subsection{Identification of $S_{\ga}$ with the Mars-Springer slice}
Let $B$ be the Borel subgroup of upper triangular matrices in $GL_n$, and let $T$ be the maximal torus
of diagonal matrices of $GL_n$.
Recall from Section~2 that in order to define the slice, we must 
choose a point $x_{\ga} = xK/K \in Q_{\ga}$ such that $x \theta(x)^{-1} \in N$, where $N$ is 
the normalizer of $T$,
and where $\theta$ is the involution (\ref{eqn:theinvolnov29}).  Here, we choose $x$ to be the 
the origin of $S_{\ga}$, meaning the matrix where 
all $\zz$-variables are set to $0$.  The fact that $x_{\ga}=xK/K$ is actually a point of $Q_{\ga}$ follows 
from \cite[Section 2.2-2]{Yamamoto-97}, while the fact that $x \theta(x)^{-1} \in N$ is easy to check.  
Denote by $X_{\ga}$ the Mars-Springer slice to
$Q_{\ga}$ at $x_{\ga}$ in $G/K$, so
\[ X_{\ga}:=(U^- \cap \psi(U^-)) \cdot x_{\ga} = ((U^- \cap \psi(U^-))x)/K. \]

\begin{Theorem}
\label{thm:bigslice}
We have an isomorphism $X_{\ga} \cong S_{\ga}$.
\end{Theorem}
\begin{proof}
Recall from Section \ref{sec:attractive-slices} that the map $\psi: G \to G$ is defined by $c_{\eta} \circ \theta$,
where $\eta = x \theta(x)^{-1} \pmod{T}$.

Given a clan $\gamma$, let $\mathcal{I}(\alpha)$ denote the {\bf underlying involution} for $\alpha$;
this is the permutation in $S_n$ such that $\mathcal{I}(\alpha)(j)=j$ if $\alpha_j=\pm$, while, if $(i<j)$
is a matching of $\alpha$, $\mathcal{I}(\alpha)(i)=j$ and $\mathcal{I}(\alpha)(j)=i$.  The map
$\mathcal{I}$ is the concrete realization for $(G,K)=(GL_n,GL_p\times GL_q)$ of the map defined by
R.~W.~Richardson--T.~A.~Springer~\cite{Richardson-Springer-90} associating a twisted involution
(which is an involution in this case) to every orbit.

\begin{Claim}
\label{claim:etadec14}
As an element of $\caW=S_n$, we have $\eta=\mathcal{I}(\alpha)$.
\end{Claim}
\begin{proof}
The matrix computation $x\theta(x)^{-1} = x I_{p,q} x^{-1} I_{p,q}$ is straightforward, and results in a monomial matrix
representing the involution described above.
\end{proof}

\begin{Claim}
\label{claim:thegroupdec14}
The group $U^- \cap \psi(U^-)$ consists of the unipotent matrices having $1$'s on the diagonal, $0$'s above the diagonal,
$0$'s below the diagonal in all positions $(i,j)$ with $i>j$ such that $\eta(i) < \eta(j)$, and arbitrary entries in all
other positions below the diagonal.
\end{Claim}
\begin{proof}
$U^-$ is the set of of lower triangular unipotent $n\times n$
matrices.   $U^-$ is $\theta$-stable, so $\theta$ maps it isomorphically
to itself, meaning that $\psi(U^-)=c_{\eta}(U^-)$.  By Claim~\ref{claim:etadec14}, $\eta=\eta^{-1}$. Thus
conjugation by $\eta$ sends
the affine coordinate at position $(i,j)$ ($i>j$) to position $(\eta(i),\eta(j))$.  The claim follows.
\end{proof}

Denote by $S_{\alpha}^\prime$ the set of $n\times n$ matrices that satisfy (O.1), (O.2), (O.3), (Z.2), (Z.4), (Z.6),
(Z.8), (Z.9), and the following condition:

\begin{itemize}
\item[(Z.10)] If $(i<j)$ is any matching of $\ga$, and $i < k < j$, then 
\begin{itemize}
  \item[$\bullet$] If $k$ is the left end of a matching $(k < \ell)$ with $j < \ell$, then the entries at positions $(k,w_{\ga}(i))$ and 
$(k,w_{\ga}(j))$ are arbitrary, but equal.
  \item[$\bullet$] Otherwise, the entries at positions $(k,w_{\ga}(i))$ and $(k,w_{\ga}(j))$ are both zero.
\end{itemize}
\end{itemize}

\begin{Claim}
\label{claim:theyequaldec14}
$S_{\alpha}^\prime=(U^- \cap \psi(U^-))x$.
\end{Claim}
\begin{proof}
By a direct but tedious matrix computation using Claim~\ref{claim:thegroupdec14},
one checks that $(U^- \cap \psi(U^-))x \subseteq S_{\ga}^\prime$.  The space $(U^- \cap \psi(U^-))x$ is
evidently an affine space of dimension $\binom{n}{2}-\ell(\eta)$, where $\ell( \cdot )$ denotes Coxeter length.
The space $S_{\ga}^\prime$ is also visibly an affine space, and a straightforward combinatorial argument
shows that its dimension is also $\binom{n}{2}-\ell(\eta)$.  Thus the aforementioned containment is in fact an equality.
\end{proof}

\begin{Claim}
\label{claim:containeddec15}
$S_{\alpha}\subseteq S_{\alpha}^{\prime}$.
\end{Claim}
\begin{proof}
This follows from the fact that (Z.3) and (Z.7) together imply (Z.10).
\end{proof}

The right action of $K$ on $G$ is by column operations which separately
act on the first $p$ columns and the last $q$ columns of any matrix in $G$.
Given a matrix $M^\prime \in S_{\ga}^\prime$, we now describe an algorithm for column
reducing $M^\prime$ via this $K$-action to another matrix $M$ of a particular form.

\begin{Algorithm}
\label{alg:col-reduce}
Begin with the first $p$ columns of $M^\prime$.  For each row $i=1,2,\ldots,n$, determine whether there is a
pivot in row $i$.  If not, move to the next row.  If so, then this pivot is at $(i,w_{\ga}(i))$.  Now start
with column $j=1$, and move right toward column $j=w_{\ga}(i)-1$.  For each such $j$, if entry $(i,j)$ is zero, 
move to the next column.  Otherwise, entry $(i,j)$ is $\lambda_{i,j} \neq 0$. Replace column $C_j$ with
$C_j - \lambda_{i,j} \cdot C_i$.  
Repeat the same procedure using the last
$q$ columns.
\end{Algorithm}

Let $M$ be the output of Algorithm~\ref{alg:col-reduce}.
\begin{Claim}
\label{Claim:leftward-operations}
The first $p$ columns and the last $q$
columns of $M$ are in reduced column echelon form, meaning that all entries in the same row as a pivot (except for the pivot) are zero.
$M$ is the unique matrix in the $K$-orbit of $M^\prime$ with this property.
\end{Claim}
\begin{proof}
Consider first the leftmost $p$ columns.  Algorithm~\ref{alg:col-reduce} brings to zero any entry in the same
row as a pivot at $(i,w_{\ga}(i))$ and to the left of 
$(i,w_{\ga}(i))$.  Since $M^\prime$ satisfies (Z.2), any entry to the right of a pivot is
already zero: there is a pivot in every column, and by our
choice of $w_{\ga}$, any pivot right of column $w_{\ga}(i)$ occurs in a row $j > i$.

For the rightmost $q$ columns, one argues identically, using the fact that $M^\prime$ satisfies (Z.6).

Uniqueness of $M$ follows from the uniqueness of reduced column echelon form (for each of the two submatrices of $M$).
\end{proof}

\begin{Claim}
\label{claim:bigonedec14}
The matrix $M$ is an element of $S_{\ga}$.
\end{Claim}
\begin{proof}
By Claim~\ref{Claim:leftward-operations}, $M$ satisfies (Z.1) and (Z.5).  We show that each step 
performed by Algorithm~\ref{alg:col-reduce} preserves properties (O.1), (O.2),
(O.3), (Z.2), (Z.4), (Z.6), and (Z.8), all of which are satisfied by $M^\prime$.  Hence $M$ also possesses
these properties.  We also show that once
Algorithm~\ref{alg:col-reduce} is complete, the resulting matrix also satisfies (Z.3), (Z.7), and (Z.9).
This will show that $M$ possesses all the properties (O.1)-(O.3) and (Z.1)-(Z.9).

(O.1):  Suppose that $c_i=+$. There is a pivot $1$ at $(i,w_{\ga}(i))$.  Any step 
of Algorithm~\ref{alg:col-reduce} affecting column $C_{w_{\ga}(i)}$ is performed to make zero
an entry $\lambda\neq 0$ at $(j,w_{\ga}(i))$ for $j>i$ (using the pivot in row $j$). Such a step replaces column
$C_{w_{\ga}(i)}$ by $C_{w_{\ga}(i)} - \lambda \cdot C_{w_{\ga}(j)}$.  Since $j>i$, by (Z.2), the entry
at $(i,w_{\ga}(j))$ is zero, so the $1$ at $(i,w_{\ga}(i))$ remains a $1$ after this column operation.

When $c_i=-$, one argues identically using the last $q$ columns and (Z.6).  Thus (O.1) is preserved by
all steps of our algorithm, as desired. 

(O.2):  The argument is identical to that for (O.1).

(O.3):  Suppose $(i<j)$ is a matching. First consider the $-1$ at $(j,w_{\ga}(i))$.  Since entries to the right of this
$-1$ are zero by (Z.4), the $-1$ is clearly unchanged by steps of the algorithm.  Similarly, the $1$ at $(j,w_{\ga}(j))$ is
unchanged, by (Z.8).  

(Z.2):  Given a zero entry above a pivot $1$, any entry to its right is also zero by (Z.2) since such an entry is
above a different pivot $1$.  Thus any step of the algorithm leaves the first of these two zero entries unchanged, as needed.

(Z.3):  Suppose $(i<j)$ is a matching. There are a $1$ and a $-1$ in positions $(i,w_{\ga}(i))$ and
$(j,w_{\ga}(i))$, respectively.  Let $i<k<j$.  We must check that $M$ is zero at $(k,w_{\ga}(i))$.  There are two cases based
on (Z.10):

\noindent
{\sf Case 1:} (The  $(k,w_{\ga}(i))$ entry of $M^\prime$ is arbitrary.)
This case occurs if and only if $\ga_k$ is the left end of a matching $(k < \ell)$ such that
$i < k < j < \ell$.  Thus there is a pivot $1$ in row $k$ at position $(k,w_{\ga}(k))$.  Since $M$ satisfies (Z.1),
entry $(k,w_{\ga}(i))$ of $M$ must be zero.

\noindent
{\sf Case 2:} (The $(k,w_{\ga}(i))$ entry of $M^\prime$ is required to be zero.)
Consider a step of Algorithm~\ref{alg:col-reduce} altering column $C_{w_{\ga}(i)}$. It does so by replacing
$C_{w_{\ga}(i)}$ with
$C_{w_{\ga}(i)}-\lambda \cdot C_{w_{\ga}(\ell)}$ for $\ell$ such that $\ell > i$ and $w_{\ga}(\ell) > w_{\ga}(i)$.
This step is done if and only if there is a pivot $1$ at $(\ell,w_{\ga}(\ell))$ and
a nonzero entry $\lambda$ in position $(\ell,w_{\ga}(i))$.  We examine the effect of this step
on entry $(k,w_{\ga}(i))$.

If $\ell > k$, then entry $(k,w_{\ga}(\ell))$ is zero by (Z.2), so entry $(k,w_{\ga}(i))$
remains zero after the step in question.  So suppose that $i < \ell < k$.  If entry $(\ell,w_{\ga}(i))$ is
$\lambda \neq 0$, then by (Z.10) this must be because
$\ga_{\ell}$ is the left end of a matching $(\ell < r)$ such that $i<\ell<j<r$.
If the entry at $(k,w_{\ga}(\ell))$ is zero, then the entry at $(k,w_{\ga}(i))$ remains zero after our
step, and we are done.  So suppose it is not.
Then by (Z.10) again, $\ga_k$ is the left end of a matching
$(k<s)$ with $\ell<k<r<s$.  But then $i<k<j<s$, contradicting our assumption that (Z.10) required
entry $(k,w_{\ga}(i))$ of $M^\prime$ to be zero.  Thus this situation cannot occur.

(Z.4):  Suppose that $M^\prime$ has a $-1$ at position $(i,j)$ with $1 \leq j \leq p$.  Since $M^\prime$ satisfies
(Z.4), any entry in a position $(i,k)$ with $j < k \leq p$ is zero at the outset.  Moreover, again by (Z.4), any entry
in a position $(i,\ell)$ with $k < \ell \leq p$ is also zero at the outset.  Thus all steps of
Algorithm~\ref{alg:col-reduce} leave the $0$ in position
$(i,k)$ unchanged.  

(Z.6):  We argue as with (Z.2) above, but using the rightmost $q$ columns.

(Z.7):  We argue as with (Z.3).

(Z.8):  We argue as with (Z.4).

(Z.9):  Suppose $\alpha$ has a  $1212$ pattern
consisting of matchings $(i<j)$ and $(k < \ell)$ with
$i < k < j < \ell$.  Then by (Z.10), in positions $(k,w_{\ga}(i))$ and $(k,w_{\ga}(j))$ of $M^\prime$, we have arbitrary but
equal entries; assume that both entries are equal to $\lambda_1$.  Moreover, by (Z.9), in positions $(\ell,w_{\ga}(i))$
and $(\ell,w_{\ga}(j))$, we have arbitrary entries that are negatives of one another.  Suppose that the entry
at $(\ell,w_{\ga}(i))$ is $\lambda_2$, and that the entry at $(\ell,w_{\ga}(j))$ is $-\lambda_2$. Consider
two steps performed by Algorithm~\ref{alg:col-reduce}:  
\begin{itemize}
\item[(i)] Replace column
$C_{w_{\ga}(i)}$ by
$C_{w_{\ga}(i)} - \lambda_1 \cdot C_{w_{\ga}(k)}$ to make zero the entry at $(k,w_{\ga}(i))$.
\item[(ii)] Replace column $C_{w_{\ga}(j)}$ by $C_{w_{\ga}(j)} - \lambda_1 \cdot C_{w_{\ga}(\ell)}$ to make zero the
entry at $(k,w_{\ga}(j))$.  
\end{itemize}
Step (i) replaces the $\lambda_2$ at $(\ell,w_{\ga}(i))$ by $\lambda_2 + \lambda_1$ (since there is a $-1$ in position $(\ell,w_{\ga}(k))$).
Meanwhile, step (ii) replaces the $-\lambda_2$ at $(\ell,w_{\ga}(j))$ by $-\lambda_2 - \lambda_1$ (since there is a $1$ in
position $(\ell,w_{\ga}(\ell))$).  Thus after (i) and (ii),
(Z.9) holds for the
$1212$ pattern $i<k<j<\ell$.

To conclude (Z.9) holds, it remains to show:

\begin{Subclaim}
Steps (i) and (ii) are the only steps of Algorithm~\ref{alg:col-reduce}
that can change either of the entries $(\ell,w_{\ga}(i))$ or $(\ell,w_{\ga}(j))$.
\end{Subclaim}
\begin{proof}  
Suppose there is a pivot $1$ (among the first $p$ columns) at $(r,w_{\ga}(r))$, and 
suppose that there is a nonzero entry at $(r,w_{\ga}(i))$.  There are two 
possibilities.  First, suppose $w_{\ga}(r) > w_{\ga}(k)$.  In this case, $(\ell,w_{\ga}(r))$ must be zero by (Z.4).  Thus a column operation involving columns 
$w_{\ga}(i)$ and $w_{\ga}(r)$ leaves $(\ell,w_{\ga}(i))$ unchanged, as desired.  Otherwise, $w_{\ga}(i) < w_{\ga}(r) < w_{\ga}(k)$.  Then $i < r < j$, so by (Z.10) 
the entry at $(r,w_{\ga}(i))$ can be nonzero only if there is a $1212$ 
pattern formed by matchings $(i<j)$ and $(r<s)$ with $i<r<j<s$. However, this situation is covered by
the argument preceding this subclaim.  A similar analysis of the last $q$ columns (using either (Z.8) or (Z.10)) 
shows that there are no steps of our algorithm that alter position
$(\ell,w_{\ga}(j))$. \end{proof}

This completes the proof of Claim~\ref{claim:bigonedec14}.
\end{proof}

By Claim~\ref{claim:bigonedec14}, we can now define a map
\[ \phi: X_{\ga} \rightarrow S_{\ga} \]
which sends the point $M^\prime \cdot K$ of $X_{\ga}$ to the matrix $M$, retaining the notation above.
Its inverse, as a set map, is the restriction of the natural
projection $\pi: G \to G/K$ to $S_{\ga}$.  (Note that this restriction does take values in $X_{\ga} = S_{\ga}^\prime/K$ by
Claim~\ref{claim:containeddec15}.)  Abusing notation, we denote this restriction simply by $\pi$.  Then letting $M^\prime$
and $M$ be as above, it is indeed clear that
\[ \pi(\phi(M^\prime \cdot K)) = \pi(\phi(M \cdot K)) = \pi(M) = M \cdot K = M^\prime \cdot K \]
and that
\[ \phi(\pi(M)) = \phi(M \cdot K) = M, \]
so that $\pi$ and $\phi$ are inverse as set maps.

It remains to check that $\phi$ and $\pi$ are actually morphisms of algebraic varieties.
Clearly $\pi$ is a morphism, since it is a restriction of the natural map $G \rightarrow G/K$.
The map $\phi$, on the other hand, is the restriction of an algebraic section of $\pi$ to $X_{\ga}$.  Indeed, it takes a
splitting $V_1 \oplus V_2$ and produces a matrix $M$ whose first $p$ columns are a basis for $V_1$ and whose last $q$
columns are a basis for $V_2$.  These two bases are of the form described in Subsection~\ref{sec:G/K-patches} when we
take $S$ to be the positions of the $+$'s and left endpoints of $\ga$ and $T$ to be the position of the $-$'s and left
endpoints of $\ga$.  Thus $X_{\ga} \subseteq (U_S \times U_T) \cap G/K$, and indeed the map $\phi$ is the
restriction to $X_{\ga}$ of the section $\phi_{S,T}: (U_S \times U_T) \cap G/K \to G$ described in
Subsection~\ref{sec:G/K-patches}.  This completes the proof.
\end{proof}

\section{The Mars-Springer variety and ideal}
\label{sec:ms-variety}
For $\gamma \in {\tt Clans}_{p,q}$, define the following collections of nonnegative integers \cite{Wyser-15}:
\begin{enumerate}
	\item $\gamma(i;+) :=$ number of matchings and $+$ signs among the first $i$ positions of $\gamma$ (for $i=1,\hdots,n$);
	\item $\gamma(i;-) := $ number of matchings and $-$ signs among the first $i$ positions of $\gamma$ (for $i=1,\hdots,n$); and
	\item $\gamma(i;j) := $ number of matchings whose left endpoint is weakly left of $i$ and whose right endpoint is strictly
	right of $j$ (for $1 \leq i < j \leq n$).
\end{enumerate}

Given any invertible $n\times n$ matrix $M$ and any pair of indices $i<j$, we define an auxiliary $n\times (i+j)$ matrix
$M^{[i;j]}$ as follows.  The first $i$ columns of $M^{[i;j]}$ are formed by taking the first $i$ columns of $M^{-1}$ and
zeroing out all entries in the last $q$ rows.  The last $j$ columns of $M^{[i;j]}$ are the first $j$ columns of
$M^{-1}$, unaltered.

\begin{Example}
Consider $\ga=1+--+1\in {\tt Clans}_{3,3}$.  Then:
\[ M_{1+--+1}(\zz) =
 \begin{pmatrix}
    1 & 0 & 0 & 1 & 0 & 0 \\
    0 & 1 & 0 & 0 & 0 & 0 \\
    0 & z_{3,2} & 0 & 0 & 1 & 0 \\
    0 & z_{4,2} & 0 & 0 & 0 & 1 \\
    0 & 0 & 1 & 0 & z_{5,5} & z_{5,6} \\
    -1 & 0 & 0 & 1 & 0 & 0
 \end{pmatrix}
,
\]
as the reader can verify.  The inverse of this matrix is
\[
M_{1+--+1}(\zz)^{-1} =
 \begin{pmatrix}
    \frac{1}{2} & 0 & 0 & 0 & 0 & -\frac{1}{2} \\
    0 & 1 & 0 & 0 & 0 & 0 \\
    0 & z_{3,2}z_{5,5}+z_{4,2}z_{5,6} & -z_{5,5} & -z_{5,6} & \frac{1}{2} & 0 \\
    \frac{1}{2} & 0 & 0 & 0 & 0 & \frac{1}{2} \\
    0 & -z_{3,2} & 1 & 0 & 0 & 0 \\
    0 & -z_{4,2} & 0 & 1 & 0 & 0
 \end{pmatrix}
.
\]

To obtain the auxiliary matrix $M_{1+--+1}(\zz)^{[2;4]}$, one
takes the first $2$ columns of $M_{1+--+1}(\zz)^{-1}$, zeroes 
out their last $q=3$ rows, and then concatenates to this the first $4$ columns of $M_{1+--+1}(\zz)^{-1}$, unaltered.  The
result is
\[
M_{1+--+1}(\zz)^{[2;4]}=
 \begin{pmatrix}
    \frac{1}{2} & 0 & \frac{1}{2} & 0 & 0 & 0 \\
    0 & 1 & 0 & 1 & 0 & 0 \\
    0 & z_{3,2}z_{5,5}+z_{4,2}z_{5,6} & 0 & z_{3,2}z_{5,5}+z_{4,2}z_{5,6} & -z_{5,5} & -z_{5,6} \\
    0 & 0 & \frac{1}{2} & 0 & 0 & 0 \\
    0 & 0 & 0 & -z_{3,2} & 1 & 0 \\
    0 & 0 & 0 & -z_{4,2} & 0 & 1
 \end{pmatrix}
.
\]
\qed
\end{Example}

Now define the \textbf{Mars-Springer variety} $\caN_{\gamma,\ga}$ to be the reduced subscheme of $S_{\ga}$ whose points
are the matrices $M$ satisfying all of the following rank conditions:

\begin{enumerate}
 \item[(R.1)] The rank of the southwest $(n-i) \times p$ submatrix of $M$ is at most $p-\gamma(i;+)$ (for $i=1,\hdots,n$);
 \item[(R.2)] The rank of the southeast $(n-i) \times q$ submatrix of $M$ is at most $q-\gamma(i;-)$ (for $i=1,\hdots,n$); and
 \item[(R.3)] The rank of $M^{[i;j]}$ is at most $j+\gamma(i;j)$ (for all $1 \leq i < j \leq n$).
\end{enumerate}

Note that $\caN_{\gamma,\ga}$ is clearly defined \textit{set-theoretically} by a certain collection of determinants.
Indeed, define the \textbf{Mars-Springer ideal} $\caI_{\gamma,\alpha}$ to be the ideal of $\C[\zz]$ generated by:

\begin{enumerate}
 \item All minors of the southwest $(n-i) \times p$ submatrix of $M_{\ga}(\zz)$ of size $p-\gamma(i;+)+1$ (for $i=1,\hdots,n$);
 \item All minors of the southeast $(n-i) \times q$ submatrix of $M_{\ga}(\zz)$ of size $q-\gamma(i;-)+1$ (for $i=1,\hdots,n$); and
 \item All minors of $M_{\ga}(\zz)^{[i;j]}$ of size $j+\gamma(i;j)+1$ (for all $1 \leq i < j \leq n$).
\end{enumerate}

(By Lemma \ref{lemma:isinvertible},
the minors of (3) above are in fact polynomials.)

\begin{Example}
 Consider the ideal $\caI_{123231,1+--+1}$.  One can check that the conditions defining $\caN_{123231,1+--+1}$ are all
 actually non-conditions with the exception of the condition that the rank of $M^{[2;4]}$ be at most $5$.  Thus the ideal
 $\caI_{123231,1+--+1}$ is generated by the determinant of the $6 \times 6$ matrix $M_{1+--+1}(\zz)^{[2;4]}$ computed above.
 One checks that this determinant equals (up to a scalar multiple)  $z_{3,2}z_{5,5}+z_{4,2}z_{5,6}$.  Thus
 \[\caI_{123231,1+--+1} = \langle z_{3,2}z_{5,5}+z_{4,2}z_{5,6}\rangle,\]
which defines a singular quadric hypersurface.
 \qed
\end{Example}

\begin{Theorem}\label{thm:MS-variety-is-slice}
 ${\mathcal N}_{\gamma,\alpha} \cong W_{\gamma} \cap X_{\ga}$, the Mars-Springer slice to $Q_{\ga}$ at $x_{\ga}$ in
 $W_{\gamma}$.
\end{Theorem}
\begin{proof}
Recall that points of $G/K$ are splittings ${\mathcal S}=V_1\oplus V_2$
of ${\mathbb C}^{p+q}$ where $\dim(V_1)=p$ and $\dim(V_2)=q$.  The $B$-orbit 
closure $W_{\gamma}\subseteq G/K$  has a 
description due to the second author, which we now state. Define $\pi_{\mathcal{S}}: \C^n \rightarrow V_1$ to be the
projection onto the first summand of $\caS$ with kernel $V_2$, the second summand of $\caS$.

\begin{theorem*}[\cite{Wyser-15}]\label{thm:wyser-set-description}
Let $E_{\bullet}$ be the standard coordinate flag on $\C^n$, with $E_i$ the span of the first $i$ standard basis vectors.
Then $W_{\gamma} \subseteq G/K$ is precisely the set of splittings $\mathcal{S} = V_1 \oplus V_2$
satisfying all of the following incidence conditions with $E_{\bullet}$:
\begin{enumerate}
	\item[(C.1)] $\dim(V_1 \cap E_i) \geq \gamma(i;+)$ for $i=1,\hdots,n$;
	\item[(C.2)] $\dim(V_2 \cap E_i) \geq \gamma(i;-)$ for $i=1,\hdots,n$;
	\item[(C.3)] $\dim(\pi_{\mathcal{S}}(E_i) + E_j) \leq j + \gamma(i;j)$ for $1 \leq i < j \leq n$.
\end{enumerate}
\end{theorem*}

Let $\phi:X_{\alpha}\to S_{\alpha}$ be the isomorphism from 
Theorem~\ref{thm:bigslice}. It suffices to show that the image of
\[\phi|_{W_{\gamma}\cap X_{\alpha}}: (W_\gamma\cap X_\alpha) \to S_{\alpha}\]
is precisely ${\mathcal N}_{\gamma,\alpha}$.

The map $\phi$ bijects a splitting  
${\mathcal S}=V_1\oplus V_2\in X_{\alpha}$ to a matrix $M \in S_{\alpha}$, the span of whose first $p$ columns is $V_1$,
and the span of whose last $q$ columns is $V_2$.  Therefore it suffices to show that under this bijection,
${\mathcal S}$ satisfies  (C.1)-(C.3) if and only if $M=\phi(\mathcal{S})$ satisfies (R.1)-(R.3).

The equivalence of (R.1) to (C.1) and of (R.2) to (C.2) is immediate.  

To see the equivalence of (C.3) to (R.3) we argue as follows. (C.3) states that
\[\dim(\pi_{\mathcal S}(E_i)+E_j) \leq j+\gamma(i;j).\] 
If $M\in S_{\ga}$ represents a splitting $\mathcal{S}$, then the matrix for $\pi_{\mathcal{S}}$ can be written as $MPM^{-1}$, where $P$ is the diagonal matrix with $p$-many $1$'s on the diagonal followed by $q$-many $0$'s.  
Hence $\pi_{\mathcal S}(E_i)+E_j$ is the span of the first $i$ columns of $MPM^{-1}$ and the
first $j$ columns of the identity matrix $I$.  Therefore, 
\begin{equation}
\label{eqn:anequivdec15}
\mathcal{S} \text{\ satisfies (C.3)} \iff
\text{rank}(M^{\overline{[i;j]}})\leq j+\gamma(i;j),
\end{equation}
 where $M^{\overline{[i;j]}}$ is the $n\times (i+j)$ matrix whose first $i$ columns are the first $i$ columns of $MPM^{-1}$ and whose last $j$ columns are the first $j$ columns of $I$. 

Since $M^{-1}$ is invertible,
\begin{equation}
\label{eqn:anequalrankdec15} 
\text{rank}(M^{[i;j]})=\text{rank}(M^{-1}M^{\overline{[i;j]}})=
\text{rank}(M^{\overline{[i;j]}}).
\end{equation}
Combining (\ref{eqn:anequivdec15}) and (\ref{eqn:anequalrankdec15}) we
conclude
\begin{eqnarray}\nonumber
\mathcal{S} \text{\ satisfies (C.3)} & \iff & 
\text{rank}(M^{[i;j]})\leq j+\gamma(i;j) \text{\ for all $i,j$ with $1\leq i<j\leq n$}\\ \nonumber
& \iff & M \text{\ satisfies (R.3)},\nonumber
\end{eqnarray}
as desired.
\end{proof}

\begin{Corollary}\label{cor:MS-variety-irreducible}
 For any clans $\ga \leq \gamma$, the Mars-Springer variety $\caN_{\gamma,\ga}$ is irreducible.
\end{Corollary}
\begin{proof}
 This is immediate by Theorem~\ref{thm:MS-variety-is-slice} combined with Corollary~\ref{cor:irreducible}.
\end{proof}

The key to our proof of Theorem~\ref{thm:main} is:
\begin{Theorem}[Isomorphism Theorem]
\label{thm:iso}
If $[\beta,\theta]$ interval pattern contains $[\alpha,\gamma]$, then
\[ \caN_{\theta,\beta} \cong \caN_{\gamma,\alpha}. \]
\end{Theorem}

This will be proved in Section~6.

\section{Combinatorics of interval embeddings}

Denote Bruhat order on clans by $<$.
Let $\gamma \lessdot \tau$ denote a covering relation, so $\gamma\lessdot\tau$ if
$\gamma < \tau$ and there is no clan $\beta$ such that
$\gamma < \beta < \tau$.  Equivalently, $\gamma\lessdot\tau$ if $\gamma < \tau$ and $\ell(\tau) - \ell(\gamma) = 1$.

We need the following result of the second author describing the covering relations $\lessdot$:

\begin{Theorem}[\cite{Wyser-15}]
\label{thm:wyserclosure}
Let $\gamma,\tau\in{\tt Clans}_{p,q}$ such that $\gamma < \tau$. Then there exists $\gamma'\in {\tt Clans}_{p,q}$
such that $\gamma < \gamma' \leq \tau$, where $\gamma'$ is obtained from $\gamma$ by one of the following
operations on patterns in $\gamma$:
\begin{itemize}
\item[(T.1)]  $+-\mapsto 11$
\item[(T.2)] $-+\mapsto 11$
\item[(T.3)] $11+\mapsto 1+1$
\item[(T.4)] $11-\mapsto 1-1$
\item[(T.5)] $+11\mapsto 1+1$
\item[(T.6)] $-11\mapsto 1-1$
\item[(T.7)] $1122\mapsto 1212$
\item[(T.8)] $1122\mapsto 1+-1$
\item[(T.9)] $1122\mapsto 1-+1$
\item[(T.10)] $1212\mapsto 1221$
\end{itemize}
\end{Theorem}

Call any ordered pair $\gamma\mapsto \gamma'$ of clans obtained by one of the operations (T.1)--(T.10) a {\bf transposition}.
Thus

\begin{Corollary}[\cite{Wyser-15}]
$\gamma \lessdot \gamma'$ if and only if $\gamma\mapsto \gamma'$ is a transposition such that $\ell(\gamma')=\ell(\gamma)+1$.
\end{Corollary}

For any clan $\gamma$, denote by $M_{\gamma}$ the set of matchings of $\gamma$.  Call a matching $(a<b) \in M_{\gamma}$ with $a<i<b<j$ {\bf incoming} to $(i<j) \in M_{\gamma}$.  For $(i<j) \in M_{\gamma}$, let
\[I(i,j,\gamma):=\#\{(a<b) \in M_{\gamma} \mid \text{\ $(a<b)$ is incoming to $(i<j)$}\},\]
and let
\[C(i,j,\gamma):=j-i-I(i,j,\gamma).\]  
  Then \cite{Yamamoto-97}
\[ \ell(\gamma):= \displaystyle\sum_{(a<b) \in M_{\gamma}} C(a,b,\gamma).\]
We will need the following combinatorial fact:

\begin{Proposition}
\label{prop:itincreases}
Suppose $\gamma\mapsto \gamma'$ is a transposition. Then
$\ell(\gamma)<\ell(\gamma')$.
\end{Proposition}
\begin{proof}
For brevity, we prove the hardest case
(T.7). Arguments for the remaining cases of (T.1)--(T.10) are similar in nature.

Suppose $\gamma\mapsto\gamma'$ via (T.7). 
Say the $1122$ pattern occurs at positions $i<j<k<\ell$.  
In order to compare $\ell(\gamma)$ to $\ell(\gamma')$, we consider other matchings $(a<b)$ and 
compare $I(a,b,\gamma')$ with $I(a,b,\gamma)$, $I(i,k,\gamma')$ with $I(i,j,\gamma)$, and
$I(j,\ell,\gamma')$ with $I(k,\ell,\gamma)$.

There are $15$ different possible orientations of a matching $(a<b)$ relative to
$(i<j)$ and $(k<\ell)$ to consider.

\noindent
{\sf Case 1} (the configuration $a<i<j<b<k<\ell$):  Here
\[I(a,b,\gamma')=I(a,b,\gamma),\]
\[I(j,\ell,\gamma')=I(k,\ell,\gamma)+1,\]
and
\[I(i,k,\gamma')=I(i,j,\gamma)+1.\]
Indeed, $(a<b)$ was incoming to neither $(i<j)$ nor $(k<\ell)$ in $\gamma$,
but it is incoming to both $(i<k)$ and $(j<\ell)$ in $\gamma'$.  So each such matching $(a<b)$ accounts for a decrease in
length of $2$ as we transpose from $\gamma$ to $\gamma'$.

\noindent
{\sf Case 2} (the configuration $i<j<a<k<\ell<b$):  Now,
\[I(a,b,\gamma')=I(a,b,\gamma)+2.\]
Indeed, in $\gamma$ neither $(i<j)$ nor $(k<\ell)$ are incoming to $(a<b)$,
but in $\gamma'$ both $(i<k)$ and $(j<\ell)$ are incoming to it.  On the other hand, the move in this configuration does not
introduce any new incoming matchings to $(i<k)$ or $(j<\ell)$.  So each matching $(a<b)$ in this configuration also accounts for a
decrease in length of $2$ as we transpose from $\gamma$ to $\gamma'$.

\noindent
{\sf Case 3} (the configuration $i<a<j<b<k<\ell$):  When we
transpose from $\gamma$ to $\gamma'$, $(a<b)$ loses an incoming matching: $(i<j)$ is incoming to $(a<b)$ in $\gamma$, but $(i<k)$
is not incoming to $(a<b)$ in $\gamma'$. However, $(j<\ell)$ gains an incoming matching compared to $(k<\ell)$, as $(a<b)$ is not
incoming to $(k<\ell)$ in $\gamma$, but it is incoming to $(j<\ell)$ in $\gamma'$.  Thus while the transposition causes a
change in the number of matchings incoming to $(a<b)$ and in the number of matchings to which $(a<b)$ is incoming, the net effect is
zero.

By similarly arguing as in {\sf Case 3} for the remaining $12$
configurations, one sees that {\sf Case 1} and
{\sf Case 2} are the only two which result in a net change in length coming from
changes in the number of incoming matchings.

The only change in the number of incoming matchings which we have not accounted for is that after the move, $(j<\ell)$
has also added the incoming matching $(i<k)$ in $\gamma'$, whereas $(k<\ell)$ did not have the incoming matching $(i<j)$ in $\gamma$.
This causes an additional decrease by $1$ in length as we perform this transposition.

Finally, we must account for the change in total lengths of matchings caused by the transposition.  This is $(k-i) + (\ell-j) - ((j-i) + (\ell-k)) = 2(k-j)$.

Putting this all together, we conclude that
\begin{multline}\label{eq:length-diff}
 \ell(\gamma') - \ell(\gamma) = 2(k-j)-1-
 2(\#\{(a<b) \in M_{\gamma} \mid a<i<j<b<k<\ell\} \\ +  \#\{(a<b) \in M_{\gamma} \mid i<j<a<k<\ell<b\}).
\end{multline}

Now, the sum in parentheses is clearly at most $k-j-1$, which implies that the right-hand side is
at least $1$.  Thus $\ell(\gamma') > \ell(\gamma)$, as desired.
\end{proof}

Say $[\alpha,\gamma]$ {\bf merely embeds} into $[\beta,\theta]$ if the pair of intervals satisfies the definition of an interval pattern embedding,
\emph{except} possibly for the length requirement $\ell(\gamma)-\ell(\alpha)=\ell(\theta)-\ell(\beta)$.
Given any $\alpha$, we define $\Phi(\alpha)$ to be the unique clan such that $[\alpha,\gamma]$
merely embeds in $[\Phi(\alpha),\theta]$, so $\Phi(\alpha)$ agrees with $\alpha$ on the set of
embedding indices $I$ and agrees with $\theta$ on $\{1,\ldots,n\}\setminus I$.

\begin{Theorem}
\label{lemma:interval}
Assume $[\alpha,\gamma]$ merely embeds into $[\Phi(\alpha),\theta]$.
\begin{itemize}
\item[(I)] $\ell(\gamma)-\ell(\alpha)\leq \ell(\theta)-\ell(\Phi(\alpha))$.
\item[(II)] Further suppose  $[\alpha,\gamma]$ interval embeds into
$[\Phi(\alpha),\theta]$.
If $\alpha \lessdot \alpha' \leq \gamma$ then
\[\Phi(\alpha) \lessdot \Phi(\alpha') \leq \theta.\]
\end{itemize}
\end{Theorem}
\begin{proof}
(I): Pick a chain of covering transpositions 
\[\alpha=\alpha^{(0)}\mapsto \alpha^{(1)}\mapsto \cdots \mapsto \alpha^{(m)}=\gamma.\] 
This induces a chain from $\Phi(\alpha)$ to  $\theta$ by using the same transpositions
(relative to the embedding). By Proposition~\ref{prop:itincreases},
each transposition in the latter chain increases the length by at least one, from which the statement follows.

(II): There exists a covering transposition from $\alpha$ to $\alpha'$, followed by a chain of covering transpositions from
$\alpha'$ to $\gamma$. Thus $\Phi(\alpha)$ is at least related to $\Phi(\alpha')$ by the same transposition (relative to the
embedding), and the chain from
$\alpha'$ to $\gamma$ induces a chain $\Phi(\alpha')\to\cdots\to \theta$.  Hence we have 
\[\Phi(\alpha) < \Phi(\alpha') \leq \theta.\]

It remains to show 
\[\Phi(\alpha')\lessdot\Phi(\alpha).\]
By Proposition~\ref{prop:itincreases}, we have 
\[\ell(\Phi(\alpha'))-\ell(\Phi(\alpha))\geq 1.\] 
We are done unless
\[\ell(\Phi(\alpha'))-\ell(\Phi(\alpha))\geq 2,\] 
so assume this. Let 
\[k:=\ell(\gamma)-\ell(\alpha)=\ell(\theta)-\ell(\Phi(\alpha)).\]
Since $[\alpha',\gamma]$ merely embeds into  $[\Phi(\alpha'),\theta]$, we obtain a contradiction:
\begin{eqnarray}\nonumber
k-1 & = &\ell(\gamma)-\ell(\alpha') \leq \ell(\theta)-\ell(\Phi(\alpha'))\\  \nonumber
& = & \ell(\theta)-\ell(\Phi(\alpha))-(\ell(\Phi(\alpha')) - \ell(\Phi(\alpha))) \\ \nonumber
&\leq & k-2, \nonumber
\end{eqnarray}
where for the first inequality we have used (I).
\end{proof}

\begin{Corollary}
\label{cor:biggeralpha}
If $[\alpha,\gamma]\hookrightarrow [\Phi(\alpha),\theta]$
and $\alpha \lessdot \alpha' \leq \gamma$, then $[\alpha',\gamma] \hookrightarrow [\Phi(\alpha'),\theta]$.
\end{Corollary}
\begin{proof}
By Theorem~\ref{lemma:interval}(II), 
\[\Phi(\alpha) \lessdot \Phi(\alpha') \leq \theta.\]
Hence $[\alpha',\gamma]$ and $[\Phi(\alpha'),\theta]$ have the same length difference, as desired.
\end{proof}

\begin{Corollary}
If $\theta$ avoids $[\alpha',\gamma]$ then $\theta$ avoids $[\alpha,\gamma]$ for all $\alpha \leq \alpha'$.
\end{Corollary}
\begin{proof}
This is the contrapositive of Corollary~\ref{cor:biggeralpha}, combined with induction on $\ell(\alpha')-\ell(\alpha)$.
\end{proof}

For the next lemma, assume that $[\ga,\gamma]$ interval embeds into $[\beta,\theta]$,
let $I$ be the set of indices for the embedding, and let
$J=\{1,\ldots,n\}\setminus I$ be the set of indices not involved in the embedding.

\begin{Lemma}
\label{lemma:isoprooflemma1}
\begin{itemize}
\item[(I)] If $\gb_j$ (and hence $\theta_j$) is a sign for some $j\in J$, then 
\[\gb(j;+)=\theta(j;+) \text{\ and
$\gb(j;-)=\theta(j;-)$.}\]
\item[(II)] If $\gb_a$ and $\gb_{b}$ are a matched pair for some $a,b \in J$, then
\[\gb(a;b)=\theta(a;b).\]
\end{itemize}
\end{Lemma}

\begin{Example}
Let 
\[\text{$[\ga,\gamma]=[1212,1221]\subseteq {\tt Clans}_{2,2}$
and $[\gb,\theta]=[1+212,1+221]\subseteq {\tt Clans}_{3,2}$.}\]
Both are intervals of length $1$, and hence $[\alpha,\gamma]\hookrightarrow
[\beta,\theta]$. Here 
\[\text{$n=5, I=\{1,3,4,5\}$ and $J=\{2\}$.}\]  
This is an instance of Lemma~\ref{lemma:isoprooflemma1}(I) since
$\gb_2=\theta_2=+$, $\gb(2;+)=\theta(2;+)=~1$, and $\gb(2;-)=\theta(2;-)=0$.\qed
\end{Example}

\begin{Example}
Recall Example~\ref{exa:intro1}, where
\[\text{
$[\alpha,\gamma]=[+--+, 1212]$ and $[\beta,\theta]=[1+--+1, 123231]$.}\]  
Here 
\[\text{$n=6$, $I=\{2,3,4,5\}$, and $J=\{1,6\}$.}\] 
This
is an instance of Lemma~\ref{lemma:isoprooflemma1}(II) since $a=1$ and $b=6$ are matched and
$\gb(1;6)=\theta(1;6)=0$.\qed
\end{Example}

\noindent
\begin{proof}[Proof of Lemma~\ref{lemma:isoprooflemma1}]
The proof is by induction on 
\[k=\ell(\gamma)-\ell(\alpha)=\ell(\theta)-\ell(\beta).\] 
The base case, $k=0$, when
$\theta=\gb$, is trivial, so we may assume $k\geq 1$.  Suppose that 
\[\alpha \lessdot \alpha' \leq \gamma.\]
By Corollary~\ref{cor:biggeralpha}, 
\[\text{$[\alpha',\gamma] \hookrightarrow [\beta',\theta]$, where $\beta'=\Phi(\ga')$.}\]
For the statement (I), by induction, 
\[\gb'(j;\pm)=\theta(j;\pm).\]  
Next, clearly 
\[[\ga,\ga']\hookrightarrow [\gb,\gb'].\]
If we know that 
\[\gb(j;\pm)=\gb'(j;\pm),\] 
then we are done.  Hence we are reduced to the case $k=1$.  The statement (II)
similarly reduces to this situation.

Therefore we prove both claims when $k=1$.  To do this,
one must analyze each possible covering transposition
$\gb \lessdot \gb'$.  As with the proof of Proposition \ref{prop:itincreases}, for brevity we give the details
only for the hardest case (T.7). The other cases are similar.

Suppose the $1122$ pattern of $\gb$ occurs at \[i<j<k<\ell.\]   From the proof of Proposition \ref{prop:itincreases},
we recall from \eqref{eq:length-diff} that $\ell(\gb') - \ell(\gb)=1$ if and only if each position between
$j$ and $k$ is either the right endpoint of a matching $(a<b)$ with \[a<i<j<b<k<\ell\] 
or the left endpoint of a matching $(a<b)$
with \[i<j<a<k<\ell<b.\]

For (I), let $m \in J$ be such that $\gb_m=\gb'_m$ are the same sign.  If $m < j$ or
$m > k$, then 
\[\gb(m;\pm)=\gb'(m;\pm).\]  
Since $\gb \lessdot \gb'$, by the observation of the previous
paragraph, we cannot have $j<m<k$.

For (II), for $(a<b) \in M_{\gb}$, one checks that $\gb(a;b) = \gb'(a;b)$ unless $(a<b)$ is in one of the three
configurations
\[\text{$i<a<j<b<k<\ell$, \ $i<j<a<b<k<\ell$,\  or $i<j<a<k<b<\ell$.}\]  
Since $\gb \lessdot \gb'$, none of these configurations are
possible, again by the above observation.
\end{proof}

\section{Proof of Theorems~\ref{thm:iso} and~\ref{thm:main}}

Assume that $[\ga,\gamma]$ interval embeds into $[\gb,\theta]$ and that $J$ is the set of indices \textit{not} involved
in this embedding.  Let $w:=w_{\gb}$ be the permutation associated to the clan $\gb$ as described in Section \ref{sec:affine-space-s}.
\begin{Proposition}
\label{prop:projection-injective}
For any point $M \in \mathcal{N}_{\theta,\gb}$ and $j \in J$, all the entries in the $j$-th row and the $w(j)$-th column of $M$ are $0$ except for the required $\pm 1$ specified by (O.1), (O.2) or 
(O.3).
\end{Proposition}

We delay the proof of Proposition~\ref{prop:projection-injective}
until Section~6.3. Assuming it, we are ready to prove Theorems \ref{thm:iso} and \ref{thm:main} in the following two
subsections.

\subsection{Proof of Theorem~\ref{thm:iso}:} Define a map
\[\Psi: \mathcal{N}_{\theta,\gb}\rightarrow S_{\ga}\]
by deleting the rows $j$ and columns $w(j)$
for $j\in J$.  Proposition~\ref{prop:projection-injective} says that $\Psi$
is injective as a map of sets.  Hence it is injective as a map of varieties,
since it is the restriction of a linear map.

If we can show that
\begin{equation}
\label{eqn:containment}
{\rm Im}(\Psi)\subseteq \mathcal{N}_{\gamma,\ga},
\end{equation}
we are done.  Indeed, this containment is then an equality
since $\mathcal{N}_{\gamma,\ga}$ is irreducible (cf. Corollary~\ref{cor:MS-variety-irreducible}) and
since
\[\dim(\mathcal{N}_{\theta,\gb})=\ell(\theta)-\ell(\gb)=
\ell(\gamma)-\ell(\ga)
=\dim(\mathcal{N}_{\gamma,\ga}),\]
where the middle equality is by the interval embedding hypothesis.

To prove \eqref{eqn:containment}, we need to show $\Psi(m)$ satisfies (R.1), (R.2) and (R.3) whenever
$m\in\mathcal{N}_{\theta,\gb}$.  By (strong) induction on $\#J\geq 0$, we reduce to the case where $\theta$ is obtained
from $\gamma$ (or, equivalently, $\gb$ is obtained from $\ga$) by adding a single $+$, a single $-$, or a
single matching. (The base case $\#J=0$ is trivial.)

\noindent
{\sf Case 1:} ($\theta$ is obtained from $\gamma$ by adding a single $+$)  Let us suppose that $\gamma$ and $\ga$ are
$(p,q)$-clans of length $n$, so that $\theta$ and $\gb$ are $(p+1,q)$-clans of length $n+1$.  Let us further suppose
that the $+$ is added between the $(\ell-1)$th and $\ell$th characters of $\ga$ and $\gamma$.  Thus 
\[J=\{\ell\}, \text{\ and $\gb_{\ell}=\theta_{\ell}=+$.}\]
This $+$ in $\gb$ corresponds to the $\ell$th row $R$ and $w(\ell)$-th column $C$ of a matrix $m \in \caN_{\theta,\gb}$; $C$ is
among the first $p$ columns of $m$.  (It is $R$ and $C$ which are deleted by the map $\Psi$.)
The column $C$ has a single $1$, which is the only nonzero entry in $C$ and in $R$.

(R.1):  For $i=1,\hdots,n$, let $R_i^+$ denote the rank of the southwest $(n-i) \times p$ submatrix of $\Psi(m)$.  For
$j=1,\hdots,n+1$, let $\caR_j^+$ denote the rank of the southwest $(n+1-j) \times (p+1)$ submatrix of $m$.  We want
to show that 
\[R_i^+ \leq p - \gamma(i;+) \text{\ for $i=1,\hdots,n$,}\] 
knowing that 
\[\caR_j^+ \leq p+1-\theta(j;+)\text{\ for $j=1,\hdots,n+1$.}\]
There are two cases:  Either $i < \ell$, or $i \geq \ell$.

If $i < \ell$, then we know that 
\[\caR_i^+ \leq p+1-\theta(i;+).\]  Furthermore, we know that $\caR_i^+=R_i^+ + 1$ since the row
$R$ is included among the last $n+1-i$ rows of $m$ and hence contributes $1$ to the rank of the southwest $(n+1-i) \times (p+1)$
submatrix of $m$ relative to the rank of the southwest $(n-i) \times p$ submatrix of $\Psi(m)$.  Finally, we know that
\[\theta(i;+)=\gamma(i;+)\] 
since $\theta$ and $\gamma$ are the same up to position $i$.  Putting these facts together gives us
that 
\[R_i^+ = \caR_i^+-1 \leq p+1-\theta(i;+)-1 = p - \gamma(i;+),\] 
the desired conclusion.

On the other hand, if $i \geq \ell$, then we have that 
\[\caR_{i+1}^+ \leq p+1-\theta(i+1;+).\]  We also have that
\[\caR_{i+1}^+ = R_i^+\] 
since the row $R$ is now not among the last $n+1-(i+1)$ rows of $m$ and hence does not contribute $1$
to the rank of the southwest $(n+1-(i+1)) \times (p+1)$ submatrix of $m$ relative to the rank of the
southwest $(n-i) \times p$ submatrix of $\Psi(m)$.  Finally, we also clearly have that 
\[\theta(i+1; +) = \gamma(i;+)+1.\]
Putting these facts together, we see that 
\[R_i^+ = \caR_{i+1}^+ \leq p+1-\theta(i+1;+)=p+1-(\gamma(i;+)+1)=p-\gamma(i;+),\]
as desired.

(R.2):  Now, we let $R_i^-$ denote the rank of the southeast $(n-i) \times q$ submatrix of $\Psi(m)$ for $i=1,\hdots,n$,
and we let $\caR_j^-$ denote the rank of the southeast $(n+1-j) \times q$ submatrix of $m$.  We want to show that
\[R_i^- \leq q-\gamma(i;-) \text{\ for $i=1,\hdots,n$,}\] 
knowing that 
\[\caR_j^- \leq q-\theta(j;-) \text{\ for $j=1,\hdots,n+1$.}\]  
Again,
we consider the cases $i < \ell$ and $i \geq \ell$.

In this case, the map $\Psi$ simply deletes the row $R$ of $0$'s from the submatrix of $m$ formed by the last $q$ columns,
which has no effect on the rank whether this row is included among the last $n+1-i$ rows or not.  Thus for $i < \ell$, we
know that \[\caR_i^- \leq q-\theta(i;-),\] 
\[\caR_i^-=R_i^-,\] 
and
\[\gamma(i;-) = \theta(i;-),\] so
\[R_i^- \leq q-\gamma(i;-),\] as desired.

When $i \geq \ell$, we have 
\[\caR_{i+1}^- \leq q-\theta(i+1;-),\] 
\[\caR_{i+1}^- = R_i^-,\] 
and 
\[\gamma(i;-) =
\theta(i+1;-).\]  
Thus 
\[R_i^- \leq q-\gamma(i;-),\] 
again as desired.

(R.3):  We now must consider the ranks of the $(i,j)$-auxiliary matrices for all possible $i<j$.  Given $i<j$, we
denote by $\text{aux}(\Psi(m),i,j)$ the $n \times (i+j)$ auxiliary matrix formed from $\Psi(m)$, and we denote by
$R_{i,j}$ the rank of this matrix.  Similarly, we denote by $\text{aux}(m,i,j)$ the $(n+1) \times (i+j)$ auxiliary matrix
formed from $m$ and by $\caR_{i,j}$ the rank of this matrix.  We want to see that 
\[R_{i,j} \leq j + \gamma(i;j) \text{\ for
all $1 \leq i<j \leq n$,}\] 
given that 
\[\caR_{i,j} \leq j + \theta(i;j) \text{\ for all $1 \leq i < j \leq n+1$.}\]

Since $m$ is obtained from $\Psi(m)$ by adding the row $\ell$ and column $w(\ell)$,
containing a $1$ in position $(\ell,w(\ell))$ and zeros elsewhere, $m^{-1}$ is obtained from $\Psi(m)^{-1}$ by adding the
row $w(\ell)$ and the column $\ell$, containing a $1$ in position $(w(\ell),\ell)$ and zeros elsewhere.  Since the
auxiliary matrices are built from $m^{-1}$ and $\Psi(m)^{-1}$, our analysis relies primarily on this observation.

Now, there are three cases to consider:  
\[i < j < \ell,\  i < \ell \leq j, \text{\ and $\ell \leq i < j$.}\]  
When $i<j<\ell$, we
compare $\text{aux}(\Psi(m),i,j)$ to $\text{aux}(m,i,j)$.  Note that
since the lone $1$ added to $\Psi(m)^{-1}$ to form $m^{-1}$ is located in column $\ell > j$, it does not appear in the
latter auxiliary matrix; therefore, $\text{aux}(m,i,j)$ is obtained from $\text{aux}(\Psi(m),i,j)$ by simply adding a single
row of zeros.  Hence 
\[R_{i,j} = \caR_{i,j},\] and since 
\[\caR_{i,j} \leq j + \theta(i;j),\] 
and furthermore because 
\[\gamma(i;j) = \theta(i;j),\]
it follows that 
\[R_{i,j} \leq j + \gamma(i;j),\] 
as desired.

Now consider the case where $i < \ell \leq j$.  In this case, we compare $\text{aux}(\Psi(m),i,j)$ to $\text{aux}(m,i,j+1)$.
Here, the latter matrix is obtained from the former by adding a single row and column, with a lone $1$ at the
intersection of this row and column, and $0$'s elsewhere.  (Note that the $1$ does in fact appear since it is in row $w(\ell) \leq p$,
by the definition of $w$.)  Thus 
\[\caR_{i,j+1} = R_{i,j} + 1,\] 
and clearly 
\[\theta(i;j+1) =
\gamma(i;j).\]  
So since 
\[\caR_{i,j+1} \leq j+1+\theta(i;j+1),\] we have 
\[R_{i,j} = \caR_{i,j+1}-1 \leq j+1+\theta(i;j+1)-1 =
j+\gamma(i;j),\] 
as desired.

Finally, consider the case $\ell \leq i < j$.  Here we compare $\text{aux}(\Psi(m),i,j)$ to $\text{aux}(m,i+1,j+1)$.
The latter matrix is obtained from the former this time by adding one row and two columns, with the row containing two $1$'s,
one at its intersection with each of the two columns.  The row and the two columns have zeros in all other entries.  Note that
this adds $1$ (\textit{not} $2$) to the rank, so
\[\caR_{i+1,j+1}=R_{i+1,j+1}+1.\]  
Note further that 
\[\gamma(i;j) =
\theta(i+1;j+1),\] so since 
\[\caR_{i+1;j+1} \leq j+1+\theta(i+1;j+1),\] 
we have 
\[R_{i+1,j+1}=\caR_{i+1,j+1}-1 \leq
j+1+\theta(i+1;j+1)-1=j+\gamma(i;j),\] 
as required.

\noindent
{\sf Case 2:} ($\theta$ is obtained from $\gamma$ by adding a single $-$) Here the arguments are similar to
those of {\sf Case 1}; we simply interchange signs and the roles of the last $q$ columns and the first $p$ columns
in the previous argument.  We omit the details.

\noindent
{\sf Case 3:} ($\theta$ is obtained from $\gamma$ by adding a single matching) Suppose that the $(p+1,q+1)$-clan $\theta$ is obtained from $\gamma$
by adding to the $(p,q)$-clan $\gamma$ a single matching, the left endpoint of which is added between positions $\ell-1$ and
$\ell$, and the right
endpoint of which is added between positions $\ell^\prime-1$ and $\ell^\prime$.  Thus $J=\{\ell,\ell^\prime\}$.

(R.1):  As above, let $R_i^+$ denote the rank of the southwest $(n-i) \times p$ submatrix of $\Psi(m)$, and let
$\caR_j^+$ denote the rank of the southwest $(n+2-j) \times (p+1)$ submatrix of $m$.  Now the submatrix of $m$ formed by
its first $p+1$ columns differs from the submatrix of $\Psi(m)$ formed by its first $p$ columns in
that the former has two extra rows and one extra column.  The northernmost extra row (row $\ell$) contains a $1$ in column $w(\ell)$ and $0$'s
elsewhere, and the southernmost extra row (row $\ell^\prime$) contains a $-1$ in column $w(\ell)$ and $0$'s elsewhere.  Besides
the aforementioned $1$ and $-1$, column $w(\ell)$ has all of its other entries zero.

We want to see that 
\[R_i^+ \leq p-\gamma(i;+)  \text{\ for all $i$}\] 
given that 
\[\caR_j^+ \leq p+1-\theta(j;+) \text{\ for all $j$.}\]
We consider three cases:  
\[i < \ell,\ \  \ell \leq i < \ell^\prime, \text{\ \  and $\ell^\prime \leq i$.}\]  
If $i < \ell^\prime$, then
we know that 
\[\caR_i^+ \leq p+1-\theta(i;+)\] 
and that 
\[\theta(i;+)=\gamma(i;+).\]  
Also,
\[\caR_i^+=
R_i^+ + 1\] since the southwest $(n+2-i) \times (p+1)$ submatrix of $m$ contains both the aforementioned extra rows
relative to the southwest $(n-i) \times p$ submatrix of $\Psi(m)$, which causes the rank of the former to be $1$ (\textit{not} $2$)
higher than the rank of the latter (since the $1$ and $-1$ occur in the same column).  Thus
\[ R_i^+ = \caR_i^+ - 1 \leq p+1 - \theta(i;+) - 1 = p-\gamma(i;+), \]
as required.

Now, if $\ell \leq i < \ell^\prime$, then we know that 
\[\caR_{i+1}^+ \leq p+1-\theta(i+1;+)\] 
and that 
\[\theta(i+1;+) = \gamma(i;+).\]  
Note here that \[\caR_{i+1}^+ = R_i^+ + 1\] because the southwest $(n+2-(i+1)) \times (p+1)$
submatrix of $m$ still contains one of the two extra rows (namely row $\ell^\prime$) relative to the southwest
$(n-i) \times p$ submatrix of $\Psi(m)$; this causes the rank of the former to be $1$ larger than the rank of the latter.
Then
\[ R_i^+ = \caR_{i+1}^+ - 1 \leq p+1 - \theta(i+1;+) - 1 = p - \gamma(i;+). \]

Finally, if $\ell^\prime \leq i$, then we know that 
\[\caR_{i+2} \leq p+1-\theta(i+2;+)\] 
and
that 
\[\theta(i+2;+) = \gamma(i;+)+1.\]  
Here, 
\[\caR_{i+2}^+ = R_i^+\] 
because the southwest $(n+2-(i+2)) \times (p+1)$
submatrix of $m$ now does not contain either of the additional rows, so it differs from the southwest $(n-i) \times p$
submatrix of $\Psi(m)$ only in that it has an extra column of zeros.  Thus we have
\[ R_i^+ = \caR_{i+2}^+ \leq p+1-\theta(i+2;+) = p+1-(\gamma(i;+)+1) = p-\gamma(i;+), \]
as required.

(R.2):  The argument is virtually identical to that for (R.1).

(R.3):  We again must consider the ranks of the $(i,j)$-auxiliary matrices for all possible $i<j$.  As above, given $i<j$,
we denote by $\text{aux}(\Psi(m),i,j)$ the $n \times (i+j)$ auxiliary matrix formed from $\Psi(m)$, and by
$R_{i,j}$ the rank of this matrix.  We denote by $\text{aux}(m,i,j)$ the $(n+2) \times (i+j)$ auxiliary matrix
formed from $m$, and by $\caR_{i,j}$ the rank of this matrix.  We want to see that 
\[R_{i,j} \leq j + \gamma(i;j) \text{\ for
all $1 \leq i<j \leq n$}\] 
given that 
\[\caR_{i,j} \leq j + \theta(i;j)\text{\ for all $1 \leq i < j \leq n+2$.}\]

Recall that $m$ is obtained from $\Psi(m)$ by adding rows $\ell$ and $\ell^\prime$ and columns $w(\ell)$ and $w(\ell^\prime)$,
with $1$'s in positions $(\ell,w(\ell))$, $(\ell,w(\ell^\prime))$, and $(\ell^\prime,w(\ell^\prime))$, a $-1$ in position
$(\ell^\prime,w(\ell))$, and $0$'s elsewhere.  Thus $m^{-1}$ is obtained from $\Psi(m)^{-1}$ by adding rows $w(\ell)$ and
$w(\ell^\prime)$ and columns $\ell$ and $\ell^\prime$, with $1/2$'s in positions $(w(\ell),\ell)$, $(w(\ell^\prime),\ell)$,
and $(w(\ell^\prime),\ell^\prime)$, a $-1/2$ in position $(w(\ell),\ell^\prime)$, and $0$'s elsewhere.  The analysis below
relies primarily on this observation.

There are several cases to consider, depending upon the values of $i<j$ relative to $\ell$ and $\ell^\prime$.  Given $i<j$,
we define $i_{\theta}<j_{\theta}$ respectively to be the positions in $\theta$ of the $i$th and $j$th characters of $\gamma$ relative to the embedding.  So for a given $i$, there are three possibilities:
\begin{enumerate}
 \item $i_{\theta} < \ell$, in which case $i_{\theta} = i$;
 \item $\ell < i_{\theta} < \ell^\prime$, in which case $i_{\theta} = i+1$; or
 \item $\ell^\prime < i_{\theta}$, in which case $i_{\theta} = i+2$.
\end{enumerate}

Now, we consider the following cases:

($i_{\theta} < j_{\theta} < \ell$):  In this case, 
\[i_{\theta} = i, \   j_{\theta} = j, \text{\  and $\gamma(i;j) =
\theta(i_{\theta};j_{\theta})$.}\]  
Furthermore, $\text{aux}(m,i_{\theta},j_{\theta})$ differs from
$\text{aux}(\Psi(m),i,j)$ only by the addition of two rows of zeros, so 
\[\caR_{i_{\theta},j_{\theta}} = R_{i,j}.\]  
Thus
\[ R_{i,j} = \caR_{i_{\theta},j_{\theta}} \leq j_{\theta} + \theta(i_{\theta},j_{\theta}) = j + \gamma(i;j), \]
as required.

($i_{\theta} < \ell < j_{\theta} < \ell^\prime$):  Then 
\[i_{\theta} = i,\ j_{\theta}=j+1,\  \text{\ and $\gamma(i;j) =
\theta(i_{\theta},j_{\theta})$.}\]  
Here, $\text{aux}(m,i_{\theta},j_{\theta})$ differs from
$\text{aux}(\Psi(m),i,j)$ by the addition of two rows and one column.  The column contains two $1/2$'s, each of which
is contained in one of the two added rows.  All other entries of this column and these two rows are zero.  Thus
\[\caR_{i_{\theta},j_{\theta}} = R_{i,j} + 1.\]  
Then
\[ R_{i,j} = \caR_{i_{\theta},j_{\theta}} - 1 \leq j_{\theta} + \theta(i_{\theta};j_{\theta}) - 1 = j+1+\gamma(i;j)-1 =
 j+\gamma(i;j). \]

($i_{\theta} < \ell < \ell^\prime < j_{\theta}$):  Then 
\[i_{\theta} = i,\  j_{\theta} = j+2,\  \text{\ and $\gamma(i;j) =
\theta(i_{\theta},j_{\theta})$.}\]  
Here, $\text{aux}(m,i_{\theta},j_{\theta})$ differs from
$\text{aux}(\Psi(m),i,j)$ by the addition of two rows and two columns.  The first column contains two $1/2$'s,
each of which is contained in one of the two added rows.  The second column contains a $-1/2$ and a $1/2$, again with each of
these entries occurring in the two added rows.  All other entries of the added rows and columns are zero.  Thus
\[\caR_{i_{\theta},j_{\theta}} = R_{i,j} + 2.\]  
Then
\[ R_{i,j} = \caR_{i_{\theta},j_{\theta}} - 2 \leq j_{\theta} + \theta(i_{\theta};j_{\theta}) - 2 = j+2+\gamma(i;j)-2 =
 j+\gamma(i;j). \]

($\ell < i_{\theta} < j_{\theta} < \ell^\prime$):  Then 
\[i_{\theta} = i+1,\ j_{\theta} = j+1, \text{\ and $\gamma(i;j)+1 =
\theta(i_{\theta},j_{\theta})$.}\]  
Here, $\text{aux}(m,i_{\theta},j_{\theta})$ differs from
$\text{aux}(\Psi(m),i,j)$ by the addition of two rows and two columns.  The first column contains a single $1/2$
in the northmost added row.  The second column contains two $1/2$'s, with each of
these entries occurring in the two added rows.  All other entries of the added rows and columns are zero.  Thus
\[\caR_{i_{\theta},j_{\theta}} = R_{i,j} + 2.\]  Then
\[ R_{i,j} = \caR_{i_{\theta},j_{\theta}} - 2 \leq j_{\theta} + \theta(i_{\theta};j_{\theta}) - 2 = j+1+\gamma(i;j)+1-2 =
 j+\gamma(i;j). \]

($\ell < i_{\theta} < \ell^\prime < j_{\theta}$):  Then 
\[i_{\theta} = i+1,\  j_{\theta} = j+2,\text{\ and $\gamma(i;j) =
\theta(i_{\theta},j_{\theta})$.}\]  
Here, $\text{aux}(m,i_{\theta},j_{\theta})$ differs from
$\text{aux}(\Psi(m),i,j)$ by the addition of two rows and three columns; one checks that again we have
\[\caR_{i_{\theta},j_{\theta}} = R_{i,j} + 2.\]  Thus
\[ R_{i,j} = \caR_{i_{\theta},j_{\theta}} - 2 \leq j_{\theta} + \theta(i_{\theta};j_{\theta}) - 2 = j+2+\gamma(i;j)-2 =
 j+\gamma(i;j). \]

($\ell < \ell^\prime < i_{\theta} < j_{\theta}$):  Then 
\[i_{\theta} = i+2,\  j_{\theta} = j+2,\text{\  and $\gamma(i;j) =
\theta(i_{\theta},j_{\theta})$.}\]  
Here, $\text{aux}(m,i_{\theta},j_{\theta})$ differs from
$\text{aux}(\Psi(m),i,j)$ by the addition of two rows and four columns; again, it is the case that
\[\caR_{i_{\theta},j_{\theta}} = R_{i,j} + 2.\]  
Then
\[ R_{i,j} = \caR_{i_{\theta},j_{\theta}} - 2 \leq j_{\theta} + \theta(i_{\theta};j_{\theta}) - 2 = j+2+\gamma(i;j)-2 =
 j+\gamma(i;j). \]
This completes the proof of Theorem~\ref{thm:iso}, having assumed Proposition~\ref{prop:projection-injective}.\qed

\subsection{Proof of Theorem~\ref{thm:main}:} 
(I): Let $[\alpha,\gamma]\in {\mathcal C}$ be such that
$Y_{\gamma}$ is non-$\mathcal{P}$ along $\caO_{\ga}$.  Now suppose
\begin{equation}
\label{eqn:nowsupposenov23}
[\alpha,\gamma]\preceq_{\mathcal C} [\beta,\theta].
\end{equation}
Clearly, we may assume this is a covering relation, of which there are two kinds.

The first possibility is that $[\beta,\theta]$ interval pattern contains $[\alpha,\gamma]$.
By Observation~\ref{obs:orbits-equivalent}, we know that $W_{\gamma}$ is non-$\mathcal{P}$ along $Q_{\ga}$, so in
particular, $W_{\gamma}$ is non-$\mathcal{P}$
at $x_{\ga}$.  Thus $\mathcal{N}_{\gamma,\ga}$ is non-$\mathcal{P}$ at the origin, since (non-)$\mathcal{P}$ is stable
under taking slices by Lemma~\ref{lem:mildness-property-is-slice-stable}.  Then $\mathcal{N}_{\theta,\beta}$ is non-$\mathcal{P}$
at the origin as well, by Theorem~\ref{thm:iso}.  (Note that the isomorphism $\Psi$ of Theorem~\ref{thm:iso} carries the
origin on $\mathcal{N}_{\theta,\beta}$ to the origin on $\mathcal{N}_{\gamma,\ga}$.)
Thus $W_{\theta}$ is non-$\mathcal{P}$ at $x_{\gb}$, again by slice-stability of (non-)$\mathcal{P}$, and then in fact
$W_{\theta}$ is non-$\mathcal{P}$ along $Q_{\gb}$, by homogeneity.  Using Observation~\ref{obs:orbits-equivalent} once
more, we conclude that $Y_{\theta}$ is non-$\mathcal{P}$ along $\mathcal{O}_{\beta}$.

The other possibility is that $\theta=\gamma$ and $\beta\preceq \alpha$. Then 
${\mathcal O}_{\beta}\subseteq {\overline {\mathcal O}_{\alpha}}\subseteq Y_{\gamma}$, 
by the definition of the closure order. Since $Y_{\gamma}$ is non-${\mathcal P}$ along ${\mathcal O}_{\alpha}$ and
non-$\caP$ is a closed property, $Y_{\gamma}$ is non-$\caP$ on $\overline{\caO_{\ga}}$, and hence on $\caO_{\gb}$,
as desired.

(II): This is the contrapositive of (I).\qed

\subsection{Proof of Proposition~\ref{prop:projection-injective}}

Recall that $w:=w_{\gb}$ is the permutation associated to the clan $\gb$, as defined at the beginning of Section~\ref{sec:affine-space-s}.  Denote by
$\caI=\caI(\gb)$ the underlying involution of $\gb$, as defined prior to Claim~\ref{claim:etadec14}.
This means $\caI$ is the permutation that fixes $j$ if $\gb_j$ is a sign
and interchanges $i$ and $j$ if $(i<j)$ is a matching of $\gb$.

We split the proof into three main cases (A, B and C) depending on the value of $\theta_j$.

\subsubsection{Case A: ($\theta_j=+$)}  
We prove Case A via two claims.

\begin{Claim}
\label{Claim:case1firstclaim}
The entries of column $w(j)$ (except the $1$ in row $j$
 required by (O.1)) vanish.
\end{Claim}
\begin{proof}
By Lemma~\ref{lemma:isoprooflemma1}(I), $\theta(j;+)=\gb(j;+)$.  Define
\begin{equation}
\label{eqn:rsdef}
r:=p-\theta(j;+)=p-\gb(j;+).
\end{equation}
Since $\gb(n;+)=p$, it follows from the
definitions that
\begin{equation}
\label{eqn:theset123}
 r = \#\{k \mid k>j, \text{ and } \gb_k=+ \text{ or $\gb_k$ is the right end of a matching}\}. 
\end{equation}
Let $k_1<\cdots<k_r$ be the $r$ indices of the set given in
(\ref{eqn:theset123}).  Also define
\[c_i=
\begin{cases}
w(k_i) & \mbox{if $\gb_{k_i}=+$}\\
w({\mathcal I}(k_i)) & \mbox{if $\gb_k$ is the right end of a matching),}
\end{cases}\]
for $i=1,\hdots,r$.

\begin{Subclaim}
\label{Subclaim:case1}
\begin{itemize}
\item[(I)] $c_i\leq p$ for $1\leq i\leq r$.
\item[(II)] Column $c_i$ of $M$
has a fixed $\pm 1$ in row $k_i$. Besides any fixed $\pm 1$ in this row, any of the leftmost $p$ columns of $M$ strictly right of column $c_i$ has zero in row $k_i$.
\end{itemize}
\end{Subclaim}
\begin{proof}
First suppose $\gb_{k_i}=+$. Then by definition of
$w$, 
\[c_i=w(k_i) \leq p,\] 
proving (I) in this case. (II)
holds by (O.1) and (Z.1) combined.

Otherwise, $\gb_{k_i}$ is the right end of a matching. For (I), here $\caI(k_i)$ is the left end of that matching, and
\[c_i=w(\caI(k_i)) \leq p,\] 
again by $w$'s definition.  Now, by (O.2),
columns $w(\caI(k_i))$ and $w(k_i)$ of
$M$ are assigned $-1$ and $1$ respectively in row $k_i$.  (II) then holds by (Z.4).  This completes the proof of Subclaim \ref{Subclaim:case1}.
\end{proof}

Now, let ${\vec v}_i$ denote the vector consisting of the last $n-j$ entries in column $c_i$. By Subclaim~\ref{Subclaim:case1}(II),
it follows that
$\{{\vec v}_1,\ldots,{\vec v}_r\}$ is a linearly independent set.  Since $M
\in \caN_{\theta,\gb}$, by (R.1) (with $i=j$),
the southwest $(n-j) \times p$ submatrix $M^\circ$ of $M$ has rank at most
$r$ (cf.~(\ref{eqn:rsdef})).
Thus
\begin{equation}
\label{eqn:case1colspace123}
{\vec v}_1,\ldots,{\vec v}_r \text{\ is a basis of
${\rm colspace}(M^\circ)$.}
\end{equation}

Since we assume $\theta_j=+$, by (Z.2), all entries of column
$w(j)$ strictly north of row $j$ are zero. Thus it remains to show
the vector
$\vec v$ consisting of the last $n-j$ entries in column $w(j)\leq p$ is the zero vector. By
(\ref{eqn:case1colspace123}) we have
\begin{equation}
\label{eqn:case1span}
\vec v\in {\rm Span}({\vec v}_1,\ldots,{\vec v}_r).
\end{equation}
Now, for $i=1,\hdots,r$, let
\[k^\prime_i=
\begin{cases}
{\mathcal I}(k_i) & \mbox{if\ } {\mathcal I}(k_i)>j\\
k_i & \mbox{otherwise.}
\end{cases}\]
Notice that $k_i^\prime$ is either the row of a pivot, or it is a row containing a $-1$ which is located southwest of
entry $(j,w(j))$. Thus either by (Z.1) in the former case, or by (Z.4) in the latter, 
${\vec v}$ has a $0$ in row $k^\prime_i$ for all $i$. In view of (\ref{eqn:case1span}) and Subclaim~\ref{Subclaim:case1}(II), we see
that $\vec v=\vec 0$, as desired. This completes the proof of
Claim~\ref{Claim:case1firstclaim}.
\end{proof}

\begin{Claim}
The entries of row $j$ (except the $1$ in column $w(j)$ required
by (O.1)) vanish.
\end{Claim}
\begin{proof}
By (Z.1), the claim holds for the leftmost $p$
columns of $M$. Thus, it remains to check the conclusion for the rightmost $q$ columns $p+1,\hdots,n$.

Let $M^\circ$ be the southeast $n-(j-1) \times q$ submatrix of
$M$. By (R.2) for $i=j-1$, the rank of $M^\circ$ is
at most
\[q-\theta(j-1;-)=q-\theta(j;-)=q-\beta(j;-)=:r,\]
where the second equality is by Lemma \ref{lemma:isoprooflemma1}.
As in the proof of Claim~\ref{Claim:case1firstclaim}, there are $r$ positions $k_1,\hdots,k_r>j$ with $\beta_{k_i}$
being either a $-$ or the right end of a matching.  Each index corresponds  to a $1$ in row $k_i>j$
in one of the rightmost $q$ columns of $M$.  More
precisely, if $\beta_{k_i}$ is a $-$, this $1$ is the pivot of row $k_i$, and if $\beta_{k_i}$ is the right
end of a matching, the $1$ is the second $1$ in its column.
By (Z.5) and (Z.8) respectively, each of these
$1$'s have all $0$'s to their right (and in the same row).
None of these $1$'s appear in the same column,
by (O.1) and (O.3). So, these $r$ column vectors
$\vec v_1,\ldots,\vec v_r$ of $M^\circ$ are linearly
independent, and therefore
\begin{equation}
\label{eqn:isbasis345}
\vec v_1,\ldots,\vec v_r \text{\ is a basis of
${\rm colspace}(M^\circ)$.}
\end{equation}

In row $j$ of $M$, the entry in each $\vec v_i$ is zero, by (Z.6)
or (Z.7), since each such position is above either a pivot $1$ or is between two $1$'s.  So if any entry
in row $j$ among the last $q$ were nonzero,
its column in $M^\circ$ would be linearly independent of
$\{\vec v_1,\ldots,\vec v_r\}$, contradicting (\ref{eqn:isbasis345}). The claim therefore holds.
\end{proof}

\subsubsection{Case B: ($\theta_j=-$)} 
The argument is nearly identical to that of Case A; we omit the details.

\subsubsection{Case C: ($j,j^\prime\in J$ (with $j<j'$) is a matched pair)}
We start with an observation we will use repeatedly:
\begin{Claim}
\label{claim:repeatedly}
If $(k<k^\prime)$ is a matching of $\gb$, then
\[1\leq a < w(k) \implies w^{-1}(a)<k,\] 
and
\[p+1\leq a<w(k^\prime) \implies \caI(w^{-1}(a))<k.\]  
\end{Claim}
\begin{proof}
Since $k$ is the left end of a matching, $w(k)\leq p$. 
If $1\leq a < w(k) \leq p$, then by definition, $a$ is the label assigned by $w$ to either a $+$ or the left end of a
matching appearing to the left of $k$, so $w^{-1}(a)<k$. 

Similarly, since $k^\prime$ is the right end of a matching, we have
$p+1\leq w(k^\prime)\leq n$. If $p+1\leq a < w(k^\prime)$, then
$a$ is assigned by $w$ to either a $-$ or the right end $f^\prime$ of a matching $(f<f^\prime)$.  By definition of $w$, in the
former case, the ``$-$'' must appear left of $k$, so that $w^{-1}(a) < k$.  Since $\caI(w^{-1}(a)) = w^{-1}(a)$ in this
case, we have $\caI(w^{-1}(a)) < k$, as claimed.  In the latter case, where $f^\prime$ occurs at position $w^{-1}(a)$
and $f$ at position $\caI(w^{-1}(a))$, again by the definition of $w$, we must have that $f$ is left of $k$, so $\caI(w^{-1}(a)) < k$.
\end{proof}

Denote by ${\vec v}_i$ the $i$-th column vector of $M$.

Since $M$ is invertible, $\{{\vec v}_1,\ldots,{\vec v}_n\}$ is a basis of ${\mathbb C}^n$. Thus, for any $1\leq k\leq n$, we have a linear dependence relation
\begin{equation}
\label{eqn:oureqn1}
{\vec e}_k+\sum_{a=p+1}^n \lambda_{k,a} {\vec v}_a = \sum_{a=1}^p \lambda_{k,a} {\vec v}_a
\end{equation}
for some scalars $\lambda_{k,1},\ldots,\lambda_{k,n} \in {\mathbb C}$.

\begin{Claim}
\label{claim:isoproof1}
Given any $k$ such that $\gb_k=+,-,$ or the
left end of a matching, we have
\[\lambda_{k,w(b)} = \lambda_{k,w({\mathcal I}(b))}=0 \text{\ for any $b<k$.}\]
If $\gb_k$ is the right end of a matching, then
\[\lambda_{k,w(b)}=\lambda_{k,w(\caI(b))}=0 \text{\  for any $b<\caI(k)$.}\]
\end{Claim}
\noindent
\emph{Proof of Claim~\ref{claim:isoproof1}:}
Fix $k$. We prove both assertions by a common induction on $b<\text{min}\{k,\caI(k)\}$.

In the base case, when $b=1$, $\beta_1$ cannot be the right end of a matching. So we check the remaining three possibilities in turn.

If $\gb_1=+$, then the first coordinate of $\vec v_a$
is $1$ for $a=1$ and is $0$ otherwise. This follows from the definition of $w$ together with (O.1), (Z.1), and (Z.6). Since
$k>b=1$, the first row of \eqref{eqn:oureqn1} thus
says that $0=\lambda_{k,1}$. However,
\[1=w(1)=w(\caI(1))\] 
in this case, so we are done. The argument when $\gb_1=-$ is similar.

Now, suppose that $\gb_1$ is the left end of a matching. Then
the right end of this matching is at position
$\caI(1)$.  From $w$'s definition, along with (O.2), (Z.1), and (Z.5), the first
row of \eqref{eqn:oureqn1} asserts 
\[0+\lambda_{k,p+1}=\lambda_{k,1}.\]  
Meanwhile, by (O.3), (Z.4), and (Z.8), row $\caI(1)$
reads 
\[0+\lambda_{k,p+1}=-\lambda_{k,1}.\]  
Putting these together, we have 
\[\lambda_{k,w(1)}=\lambda_{k,1}=\lambda_{k,p+1}=\lambda_{k,w(\caI(1))}=0,\]
as desired.

For the inductive step, now suppose that $b>1$, and that the claims hold for indices less
than $b$.

\noindent
{\sf Case 1: ($\gb_b=+$):}
By (Z.1) in row $b$, all the entries among the first $p$ columns are $0$ except for a $1$ in column $w(b)$.  Therefore,
the $b$th row of \eqref{eqn:oureqn1} is
\begin{equation}
\label{eqn:case3asum2323}
0+\sum_{a=p+1}^n \lambda_{k,a} z_{b,a}=\lambda_{k,w(b)}.
\end{equation}
If 
\[\caI(w^{-1}(a))>b,\] 
then $z_{b,a}=0$ by (Z.6), since this entry is above the pivot in column $a$.  If on the other hand
\[\caI(w^{-1}(a))<b,\] 
then $\lambda_{k,a}=0$ by the inductive hypothesis.  Hence
$\lambda_{k,w(b)}=0$.  Since $w(\caI(b))=w(b)$ in this case, we have $\lambda_{k,w(\caI(b))} = 0$ as well.

\noindent
{\sf Case 2: ($\gb_b=-$):}  This is similar to {\sf Case 1}.

\noindent
{\sf Case 3: ($\gb_b$ is left end of a matching):} Let $b^\prime={\mathcal I}(b)$ be the right end of the
matching.
Let
\[c=w(b)  \text{\ and $c^\prime=w(b^\prime)$.}\] 
 In row $b$, \eqref{eqn:oureqn1} states that
\begin{equation}
\label{eqn:abc321}
\lambda_{k,c^\prime}=\lambda_{k,c},
\end{equation}
since
all other entries in row $b$ are $0$, by (Z.1) and (Z.5).

Now, note that if $\beta_k$ is a $+$, a $-$, or the left end of a matching, then we clearly
have that $k \neq b^\prime$.  If $\beta_k$ is the right
end of a matching, then since $b<\caI(k)$, we again have that $k\neq b^\prime$.  Thus in row $b^\prime$,
\eqref{eqn:oureqn1} states that
\begin{equation}
\label{eqn:nov10abcd}
0+\sum_{a=p+1}^{c^\prime-1} \lambda_{k,a} z_{b^\prime,a} + \lambda_{k,c^\prime} = \sum_{a=1}^{c-1} \lambda_{k,a} z_{b^\prime,a} - \lambda_{k,c},
\end{equation}
by (Z.4) and (Z.8). 

Now,
\[\caI(w^{-1}(a))<w^{-1}(c)=b\] 
for $a$ in the first sum, whereas 
\[w^{-1}(a)<w^{-1}(c)=b\] 
for $a$ in the second sum. Hence by induction, 
$\lambda_{k,a}=0$ for $1\leq a\leq c-1$ and $p+1\leq a\leq c^\prime-1$.  So
(\ref{eqn:nov10abcd}) reduces to
\[\lambda_{k,c^\prime}=-\lambda_{k,c}.\]
Combining this with (\ref{eqn:abc321}) gives  
\[\lambda_{k,c}=\lambda_{k,c^{\prime}}=0,\]
as desired.

\noindent
{\sf Case 4: ($\gb_b$ is the right end of a matching):} The left end of the matching is at position $b'=\caI(b)$. By {\sf Case 3}, 
\[\lambda_{k,w(b')}=\lambda_{k,w(\caI(b'))}=0,\] 
which is the same as 
\[\lambda_{k,w(b)}=\lambda_{k,w(\caI(b))}=0,\] 
as required.

This completes the proof.
\qed

\begin{Corollary}\label{cor:scalars-zero}
If $(k<k^\prime)$ is a matching of $\beta$, then 
\[\lambda_{k,a}=\lambda_{k^\prime,a}=0 \text{\ for $a=1,\hdots,w(k)-1$ and for $a=p+1,\hdots,w(k^\prime)-1$.}\]
\end{Corollary}
\begin{proof}
By Claim~\ref{claim:repeatedly}, for $1\leq a<w(k)$, we have $w^{-1}(a)<k$, whereas for $p+1\leq a<w(k')$,
we have $\caI(w^{-1}(a))<k$.  In either event, $\lambda_{k,a}=\lambda_{k^\prime,a}=0$ by Claim~\ref{claim:isoproof1}.
\end{proof}

\begin{Claim}
\label{claim:lambdakk+-}
If $\beta_k=+$, then $\lambda_{k,w(k)}=1$.  If $\beta_k=-$, then $\lambda_{k,w(k)}=-1$.
\end{Claim}
\noindent
\emph{Proof of Claim~\ref{claim:lambdakk+-}:}
First suppose $\beta_k=+$.  Consider row $k$ of (\ref{eqn:oureqn1}).  By (Z.1) and (O.1), this equation is
$$1+\sum_{a=p+1}^n \lambda_{k,a}z_{k,a} = \lambda_{k,w(k)}.$$
If $w^{-1}(a)<k$, then $\lambda_{k,a}=0$ by Claim~\ref{claim:isoproof1}.  On the other hand, if $w^{-1}(a)>k$, then $z_{k,a}=0$ by (Z.6) or (Z.7).  Hence the equation reduces to $1=\lambda_{k,w(k)}$, as desired.

Now, if $\beta_k=-$, then row $k$ of (\ref{eqn:oureqn1}) is
$$1+\lambda_{k,w(k)}=\sum_{a=1}^p \lambda_{k,a}z_{k,a}.$$
Again, if $w^{-1}(a)<k$, then $\lambda_{k,a}=0$ by Claim~\ref{claim:isoproof1}.  If $w^{-1}(a)>k$, then $z_{k,a}=0$ by (Z.2).
Hence the equation reduces to $1+\lambda_{k,w(k)}=0$, whence $\lambda_{k,w(k)}=-1$ as claimed.
\qed

\begin{Claim}
\label{claim:lambdakk}
Suppose $(k<k^\prime)$ is a matching of $\beta$.  Then 
\[\lambda_{k,w(k)}=1/2 \text{\ and $\lambda_{k,w(k^\prime)}=-1/2$,}\]
and
\[\lambda_{k^\prime,w(k)}=\lambda_{k^\prime,w(k^\prime)}=-1/2.\]
\end{Claim}
\noindent
\emph{Proof of Claim~\ref{claim:lambdakk}:}
First consider row $k$ of (\ref{eqn:oureqn1}).
By (O.2), (Z.1), and (Z.5),
it reads
\begin{equation}
\label{eqn:fgfg1}
1+\lambda_{k,w(k^\prime)}=\lambda_{k,w(k)}.
\end{equation}
Now consider row $k^\prime$ of (\ref{eqn:oureqn1}).
By (O.3), (Z.4), and (Z.8), it reads
\begin{equation}
\label{eqn:fgfg2}
\sum_{a=p+1}^{w(k^\prime)-1} \lambda_{k,a} z_{k^\prime,a} + \lambda_{k,w(k^\prime)} = \sum_{a=1}^{w(k)-1} \lambda_{k,a} z_{k^\prime,a}- \lambda_{k,w(k)}.
\end{equation}
Applying Corollary~\ref{cor:scalars-zero} to (\ref{eqn:fgfg2}) reduces it to
\begin{equation}
\label{eqn:fgfg3}
\lambda_{k,w(k^\prime)}=-\lambda_{k,w(k)}.
\end{equation}
Solving (\ref{eqn:fgfg1}) and (\ref{eqn:fgfg3}) simultaneously gives $\lambda_{k,w(k)}=1/2$ and $\lambda_{k,w(k^\prime)}=-1/2$,
as claimed.

Similarly, now consider \eqref{eqn:oureqn1}, but with $k$ replaced by $k^\prime$; looking at row $k$ gives
\begin{equation}
\label{eqn:fgfg5}
\lambda_{k^\prime,w(k^\prime)}=\lambda_{k^\prime,w(k)}.
\end{equation}

If we examine row $k^\prime$, we obtain
\begin{equation}
\label{eqn:fgfr4}
1+\lambda_{k^\prime,w(k^\prime)}=-\lambda_{k^\prime,w(k)}.
\end{equation}

Solving (\ref{eqn:fgfr4}) and (\ref{eqn:fgfg5}) simultaneously gives the remainder of the claim.
\qed

\begin{Claim}
\label{claim:isoproof2}
If $(k<k^\prime)$ is a matching of $\beta$, then
$$\lambda_{k,w(b)} = \lambda_{k,w({\mathcal I}(b))}=\lambda_{k^\prime,w(b)} = \lambda_{k^\prime,w({\mathcal I}(b))}=0$$
for $b=k+1,\hdots,k^{\prime}-1$.
\end{Claim}
\noindent
\emph{Proof of Claim~\ref{claim:isoproof2}:}
We only give the  proof of the claim that 
\[\lambda_{k,w(b)}=\lambda_{k,w(\caI(b))}=0.\]  
The proof that 
\[\lambda_{k^\prime,w(b)}=\lambda_{k^\prime,w(\caI(b))}=0\] 
is identical.

We induct on $b$.  For $b$ in the appropriate range, the inductive
hypothesis is that the claim holds for all $a$ satisfying $k < a < b$.

\noindent
{\sf Case 1: ($\gb_b=+$):}  In this case, row $b$ of \eqref{eqn:oureqn1} reads
\begin{equation}
\label{eqn:jeer1}
 0 + \displaystyle\sum_{a=p+1}^n \lambda_{k,a} z_{b,a} = \lambda_{k,w(b)}.
\end{equation}
Now, if $\caI(w^{-1}(a))>b$, then $z_{b,a}=0$ by (Z.6),
since this entry is above the pivot in column $a$.  Otherwise, either
\[\mathcal{I}(w^{-1}(a)) < k, \ \mathcal{I}(w^{-1}(a)) = k,\text{\ or $k < \mathcal{I}(w^{-1}(a)) < b$.}\]  
In the first case, we have
that $\lambda_{k,a} = 0$ by Claim \ref{claim:isoproof1}.  In the last case, we have that $\lambda_{k,a} = 0$ by the inductive
hypothesis.  In the middle case, we have that $z_{b,a}=0$ by (Z.7), since this entry is between the two $1$'s
in positions $(k,a)$ and $(k^\prime,a)$.  Thus the left hand side of \eqref{eqn:jeer1}
is actually zero, implying that $\lambda_{k,w(b)} =
\lambda_{k,w(\mathcal{I}(b))} = 0$, as claimed.

\noindent
{\sf Case 2: ($\gb_b = -$):} This is very similar to {\sf Case 1}; we omit the details.

\noindent
{\sf Case 3: ($\gb_b$ is the left end of a matching):}
Let the right end of this matching be $\gb_{b^\prime}$.
Let $c=w(b)$, and let $c^\prime=w(b^\prime)$.

By (Z.1) and (Z.5), in row $b$, \eqref{eqn:oureqn1} states that
\begin{equation}
\label{eqn:yyypp1}
\lambda_{k,c} = \lambda_{k,c^\prime}.
\end{equation}
By (Z.4), (Z.8) and (O.3),
in row $b^\prime$, we have
\begin{equation}
\label{eqn:c-and-c-prime}
0 + \displaystyle\sum_{a=p+1}^{c^\prime-1} \lambda_{k,a} z_{b^\prime,a} + \lambda_{k,c^\prime} =
\displaystyle\sum_{a=1}^{c-1} \lambda_{k,a} z_{b^\prime,a} - \lambda_{k,c}
\end{equation}

\noindent
{\sf Subcase 3.1: ($b^\prime < k^\prime$):}
In fact each term of the sum on the left hand side of (\ref{eqn:c-and-c-prime}) vanishes for one of three reasons:
\begin{itemize}
 \item[(L.3.1.1)] $\mathcal{I}(w^{-1}(a)) < k$: $\lambda_{k,a} = 0$ by Claim \ref{claim:isoproof1}.
 \item[(L.3.1.2)] $\mathcal{I}(w^{-1}(a)) = k$: $z_{b^\prime,a}=0$ by (Z.7), since in column
$a=w(k^\prime)$, the entry in row $b^\prime$ is located between two $1$'s in rows $k$ and $k^\prime$, since we have
assumed that $b^\prime < k^\prime$.
 \item[(L.3.1.3)] $\mathcal{I}(w^{-1}(a)) > k$: $\lambda_{k,a}=0$ by the inductive hypothesis, since $p+1 \leq a<w(b^\prime)$
 implies that $\mathcal{I}(w^{-1}(a)) < b$ by Claim~\ref{claim:repeatedly}.
\end{itemize}

Similarly, each term of the sum on the right hand side of (\ref{eqn:c-and-c-prime}) vanishes for one of these three reasons:
\begin{itemize}
 \item[(R.3.1.1)] $w^{-1}(a) < k$: $\lambda_{k,a} = 0$ by Claim \ref{claim:isoproof1}.
 \item[(R.3.1.2)] $w^{-1}(a) = k$: $z_{b^\prime,a}=0$ by (Z.3), since in column $a = w(k)$, the entry in row $b^\prime$ is located between a $1$ in row $k$
and a $-1$ in row $k^\prime$, again because we have assumed that $b^\prime < k^\prime$.
 \item[(R.3.1.3)] $w^{-1}(a) > k$: $\lambda_{k,a} = 0$ by the inductive hypothesis, since $1 \leq a < w(b)$ implies that $w^{-1}(a) < b$ by Claim~\ref{claim:repeatedly}.
\end{itemize}

Thus by (L.3.1.1)--(L.3.1.3) and (R.3.1.1)--(R.3.1.3) combined, \eqref{eqn:c-and-c-prime} reduces to
$\lambda_{k,c^\prime} = -\lambda_{k,c}$.  The simultaneous
solution of  this
with (\ref{eqn:yyypp1}) is
$\lambda_{k,c} = \lambda_{k,c^\prime} = 0$, as desired.

\noindent
{\sf Subcase 3.2: ($b^\prime > k^\prime$):}  The two matchings $(k<k^\prime)$ and $(b<b^\prime)$ form a $1212$
pattern.  We analyze each term of the sum on the left
hand side of \eqref{eqn:c-and-c-prime}:
\begin{itemize}
 \item[(L.3.2.1)] $\mathcal{I}(w^{-1}(a)) < k$:  $\lambda_{k,a} = 0$ by Claim \ref{claim:isoproof1}.
 \item[(L.3.2.2)] $\mathcal{I}(w^{-1}(a)) = k$:  We cannot conclude that $z_{b^\prime,a}=0$.  We do
know by Claim \ref{claim:lambdakk} that $\lambda_{k,a} = \lambda_{k,w(k^\prime)} = -1/2$, so all we can say is that this
particular term is equal to $-1/2z_{b^\prime,w(k^\prime)}$.
 \item[(L.3.2.3)] $\mathcal{I}(w^{-1}(a)) > k$: $\lambda_{k,a} = 0$ by the inductive hypothesis, as in (L.3.1.3).
\end{itemize}

For the right hand side of (\ref{eqn:c-and-c-prime}), we see:
\begin{itemize}
 \item[(R.3.2.1)] $w^{-1}(a) < k$: $\lambda_{k,a} = 0$ by Claim \ref{claim:isoproof1}.
 \item[(R.3.2.2)] $w^{-1}(a) = k$:  We do not know that $z_{b^\prime,a}=0$, but we do at least
know by Claim \ref{claim:lambdakk} that $\lambda_{k,a} = \lambda_{k,w(k)} = 1/2$.  Thus this term is equal to
$1/2z_{b^\prime,w(k)}$.
 \item[(R.3.2.3)] $w^{-1}(a) > k$: $\lambda_{k,a}=0$ by the inductive hypothesis, as in (R.3.1.3).
\end{itemize}

In view of the $1212$ pattern
occurring in positions $k<b<k^\prime<b^\prime$,
by (Z.9),
\[z_{b^\prime,w(k^\prime)} = -z_{b^\prime,w(k)}.\]
Thus the terms from (L.3.2.2) and (R.3.2.2) may be cancelled
in \eqref{eqn:c-and-c-prime}. Therefore
\eqref{eqn:c-and-c-prime} reduces to $\lambda_{k,c^\prime} = -\lambda_{k,c}$.  The simultaneous solution of this with
(\ref{eqn:yyypp1}) is 
\[\lambda_{k,c}=\lambda_{k,c^\prime} = 0,\] 
as required.

\noindent
{\sf Case 4: ($\gb_b$ is the right end of a matching):} If the left end of this matching occurs prior
to position $k$, we are done by Claim \ref{claim:isoproof1}.
Otherwise it occurs after position $k$, and we may apply {\sf Case 3}.  This concludes the proof of Claim \ref{claim:isoproof2}.
\qed

Let
\[ \pi_{M}:{\mathbb C}^n\to V_1\]
be the projection onto 
\[V_1:=\langle {\vec v}_1,\ldots,{\vec v}_p \rangle\]
with kernel
\[V_2:=\langle {\vec v}_{p+1},\ldots,{\vec v}_n \rangle.\]
Applying $\pi_M$ to both sides of (\ref{eqn:oureqn1}) gives
\begin{equation}
\label{eqn:oureqn2}
\pi_M({\vec e}_k)=\sum_{a=1}^p \lambda_{k,a} {\vec v}_a.
\end{equation}

\begin{Claim}
\label{claim:pinonzero}
Suppose $(k<k^\prime)$ is a matching of $\beta$.  Then
\begin{itemize}
\item[(I)] $\pi_M({\vec e}_k)$ is equal to $-1/2$ in row $k^\prime$.
\item[(II)] If $\ell$ is the left end of a matching with $\ell > k$, then
$\pi_M({{\vec e}}_{\ell})$ has entry $0$ in row $k'$.
\item[(III)] 
$\pi_M({\vec e}_k)$ is zero in row $\ell$ for any $\ell<k^\prime$ except $\ell=k$.
\end{itemize}
\end{Claim}

\noindent
\emph{Proof of Claim~\ref{claim:pinonzero}:}
(I): We have
\begin{equation}
\label{eqn:fgfg4}
\pi_M({\vec e}_k)=\sum_{a=1}^p \lambda_{k,a} {\vec v}_a=\sum_{a=w(k)}^p \lambda_{k,a}{\vec v}_{a}=\frac{1}{2}{\vec v}_{w(k)}+\sum_{a=w(k)+1}^p \lambda_{k,a} {\vec v}_{a}.
\end{equation}
The first equality is (\ref{eqn:oureqn2}).
The
second equality is by Corollary~\ref{cor:scalars-zero}.  The third equality applies Claim~\ref{claim:lambdakk}.
Now, by (O.3), $z_{k^\prime, w(k)}=-1$, and by (Z.4), $z_{k^\prime,a}=0$ for $w(k) < a \leq p$.  Therefore row $k^\prime$ of
\eqref{eqn:fgfg4} is clearly $-1/2$, as claimed.

(II):  By (Z.4) and \eqref{eqn:oureqn2}, entry $k^\prime$ of $\pi_M({\vec e}_{\ell})$ is
$\sum_{a=1}^{w(k)} \lambda_{\ell,a} z_{k^\prime,a}$.  But if $a \leq w(k)$, we have $w^{-1}(a) \leq k < \ell$ by
Claim \ref{claim:repeatedly}.  Thus by Claim \ref{claim:isoproof1}, all $\lambda_{\ell,a}$ in this sum are zero.

(III): $\pi_M(\vec e_k)$ in row $\ell$ is $\sum_{a=1}^p \lambda_{k,a}z_{\ell,a}$
 However, $\lambda_{k,a}=0$ if $w^{-1}(a)<k'$ and $w^{-1}(a)\neq k$, either by Claim \ref{claim:isoproof1}
(if $w^{-1}(a) < k$) or Claim \ref{claim:isoproof2} (if $k < w^{-1}(a) < k^\prime)$.

Now suppose $w^{-1}(a) = k$, so that $a=w(k)$.  Then there is a $1$ in position $(k,a)$ and a $-1$ in position
$(k^\prime,a)$.  If $\ell < k$, then $z_{\ell,a}=0$ by (Z.2), while if $k < \ell < k^\prime$, we have $z_{\ell,a}=0$
by (Z.3).

Finally, if $w^{-1}(a) > k^\prime$, then since $1 \leq a \leq p$ (and hence $\gb$ has either a $+$ or a left endpoint
at position $w^{-1}(a)$), there is a pivot at position $(w^{-1}(a),a)$.  Then since $\ell < k^\prime$, we have
$z_{\ell,a} = 0$ by (Z.2).  So in fact every term of the sum is zero.
\qed

Set
\[r:=\theta(j;j^\prime).\]
By Lemma~\ref{lemma:isoprooflemma1},
\[r=\gb(j;j^\prime)=\#\{(k<k^\prime) \text{\ a matching of $\beta$} \mid k<j<j^\prime<k^\prime\}.\]
Let the $r$ matchings of the above set
be
\[(k_1<k_1^\prime),(k_2<k_2^\prime),\ldots,
(k_r<k_r^\prime), \text{\ with $k_1<k_2<\cdots< k_r<j$.}\]

\begin{Claim}
\label{claim:pibasis}
\begin{itemize}
\item[(I)] The set \[\{\pi_M({\vec e}_{k_1}), \pi_M({\vec e}_{k_2}),\ldots, \pi_M({\vec e}_{k_r})\} \text{ (modulo $E_{j^\prime}$)} \] 
forms a basis of $\pi_M(E_j)+E_{j^\prime}/E_{j^\prime}$.
\item[(II)] $\pi_M({\vec e}_j)={\vec 0}$ as an element of $\pi_M(E_j)+E_{j^\prime}/E_{j^\prime}$.  In other words, for any $t>j^\prime$, $\pi_M(\vec{e_j})$ is $0$ in row
$t$.
\end{itemize}
\end{Claim}
\noindent
\emph{Proof of Claim~\ref{claim:pibasis}:}
(I):  Suppose in the quotient $\pi_M(E_j) + E_{j^\prime}/E_{j^\prime}$ that we have a dependence relation of the form
\[ \displaystyle\sum_{i=1}^r \lambda_i (\pi_M({\vec e_{k_i}}) + E_{j^\prime}) = \displaystyle\sum_{i=1}^r \lambda_i \pi_M({\vec e_{k_i}}) + E_{j^\prime} = 0 \text{ (modulo $E_{j^\prime}$).} \]
Then $\sum_{i=1}^r \lambda_i \pi_M({\vec e_{k_i}}) \in E_{j^\prime}$, meaning that for any index $\ell > j'$, the vector 
$\sum_{i=1}^r \lambda_i \pi_M({\vec e_{k_i}})$ has entry $0$ in position $\ell$.  In particular, this vector has entry $0$ in
positions $k_1^\prime,\hdots,k_r^\prime$.  Then applying parts (I) and (II) of 
Claim~\ref{claim:pinonzero}, a triangularity argument gives that all the scalars $\lambda_i$ are equal to zero.  Hence
the vectors $\pi_M({\vec e_{k_1}}),\hdots,\pi_M({\vec e_{k_r}})$ descend to a linearly independent set in the quotient
$\pi_M(E_j) + E_{j^\prime} / E_{j^\prime}$.

Then by (C.3) and the discussion following the statement of Theorem \ref{thm:wyser-set-description}, we know that
$\dim(\pi_M(E_j) + E_{j^\prime}) \leq j^\prime + \theta(j;j^\prime) = j^\prime + \gb(j;j^\prime) = j^\prime + r$.
Then $\dim(\pi_M(E_j) + E_{j^\prime}/E_{j^\prime}) \leq r$, so in fact the linearly independent set
$\{\pi_M({\vec e}_{k_1}),\hdots,\pi_M({\vec e}_{k_r})\}$ is a basis for $\pi_M(E_j)+E_{j^{\prime}}/E_{j^{\prime}}$,
as claimed.

(II):
From Part (I), we have that $\pi_M({\vec e}_j)$, viewed as an element of
$\pi_M(E_j)+E_{j^\prime}/E_{j^{\prime}}$, is a linear combination of
$\{\pi_M({\vec e}_{k_1}), \pi_M({\vec e}_{k_2}),\ldots, \pi_M({\vec e}_{k_r})\}$, or equivalently,
\[ \pi_M({\vec e_j}) + E_{j^\prime} = \displaystyle\sum_{i=1}^r \lambda_i \pi_M({\vec e}_{k_i}) + E_{j^\prime} \]
for some scalars $\lambda_1,\hdots,\lambda_r$.  We want to show that $\lambda_i=0$ for each $i$.  Since we know that
$k_i<j<k_i^\prime$ for $1\leq i \leq r$, we have by Claim~\ref{claim:pinonzero}(II) that $\pi_M({\vec e}_{j})$ has entry
$0$ in row $k_i^\prime >j^{\prime}$ for $1 \leq i \leq r$.  We also know by Claim~\ref{claim:pinonzero}(I) that
$\pi_M({\vec e}_{k_i})$ is nonzero in position $k_i^\prime$ for $1 \leq i \leq r$.  Combining these two facts
gives that each $\lambda_i=0$, and hence we have
$\pi_M(\vec e_j)=0$ in $\pi_M(E_j)+E_{j^\prime}/E_{j^\prime}$, as desired.
Since we thus have $\pi_M(\vec e_j)\in E_{j^\prime}$ the remaining assertion
also follows. \qed

We also record the following related Claim for later use:

\begin{Claim}
\label{claim:pibasisminusone}
\begin{itemize}
\item[(I)] The set 
\[\{\pi_M({\vec e}_{k_1}), \pi_M({\vec e}_{k_2}),\ldots, \pi_M({\vec e}_{k_r})\} \text{\ (modulo $E_{j^\prime-1}$)}\] forms a basis of $(\pi_M(E_{j-1})+E_{j^\prime-1})/E_{j^\prime-1}$.
\item[(II)] $\pi_M({\vec v})=0$ in row $j^\prime$ for any vector $\vec v$.
\end{itemize}
\end{Claim}
\noindent
\emph{Proof of Claim~\ref{claim:pibasisminusone}:}
(I): Since $(j<j^\prime)$ is a matching, we have 
\[\theta(j-1; j^\prime-1)=\theta(j; j^\prime)=r,\] 
since a matching $(k<k^\prime)$ satisfies 
\[k\leq j<j^\prime<k^\prime \iff k\leq j-1<j^\prime-1<k^\prime.\]  
Hence the claim follows by exactly the same argument as for Claim~\ref{claim:pibasis}(I).

(II):  By (I), it suffices to prove this for $\vec e_{k_i}$ for $i=1,\hdots,r$.  Since $j^\prime< k_i^\prime$ by definition, this follows from Claim~\ref{claim:pinonzero}(III).
\qed

Now, by (O.2),
\[z_{j,w(j)}=1 \text{ and } z_{j,w(j^\prime)}=1. \]

By (O.3),
\[ z_{j^\prime,w(j)}=-1 \text{ and } z_{j^\prime,w(j^\prime)}=1. \]

By (Z.2) and (Z.3),  
\[z_{b,w(j)}=0 \text{\ if $b<j$ or
$j<b<j^\prime$.}\] 

Finally, by (Z.6) and (Z.7),
\[z_{b,w(j^\prime)}=0 \text{\ if $b<j$ or
$j<b<j^\prime$.}\]
 
Hence the first assertion of our next claim is what remains to obtain Proposition~\ref{prop:projection-injective}'s desired
conclusion about columns $w(j)$ and $w(j^\prime)$.  The second assertion is a technical strengthening for the induction argument we give.
\begin{Claim}[Case C:  Column zeroness]
\label{claim:zbw}
$z_{b,w(j)}=z_{b,w(j^\prime)}=0$ and $\lambda_{j,w(b)}=\lambda_{j,w({\mathcal I}(b))}=0$ for all $b>j^\prime$.
\end{Claim}
\noindent
\emph{Proof of Claim~\ref{claim:zbw}:}
Our argument is by induction on $b>j'$. The inductive hypothesis is that the claims hold for all $q$ satisfying $j'<q<b$. Our inductive step is to prove that
this implies that the claims hold for $b$.
Our argument depends upon the value of $\gb_b$.

($\beta_b=+$):
 The expression for $\pi_M({\vec e}_j)$ in row $b$ (cf.~(\ref{eqn:oureqn2})) is
\begin{equation}
\label{eqn:ffqq1}
\sum_{a=1}^p \lambda_{j,a} z_{b,a}.
\end{equation}
However, 
by (O.1) and (Z.1) we know
\[z_{b,w(b)}=1 \text{\ and $z_{b,t}=0$ for $1\leq t\leq p$ and
$t\neq w(b)$.}\] 
Thus, in particular we obtain the claim's assertion that 
$z_{b,w(j)}=0$. Moreover,
(\ref{eqn:ffqq1}) reduces to 
$\lambda_{j,w(b)}$. By Claim~\ref{claim:pibasis}(II), 
all entries of $\pi_M(\vec e_j)$ beyond position $j^\prime$ are zero. Since we assume
$b>j^\prime$, we therefore conclude that 
$\lambda_{j,w(b)}=0$.  Since $\beta_b=+$, by definition we have
\[w({\mathcal I}(b))=w(b), \] 
so 
\[\lambda_{j,w({\mathcal I}(b))}=\lambda_{j,w(b)}=0,\]
proving the claim's second assertion in this case.

It remains to show that $z_{b,w(j^\prime)}=0$. For this, we use \eqref{eqn:oureqn1}, taking $k=j$, which in row $b$
gives
\begin{equation}
\label{eqn:zerosumeqn1}
0+\sum_{a=p+1}^n \lambda_{j,a} z_{b,a}=\sum_{a=1}^p \lambda_{j,a} z_{b,a}=\lambda_{j,w(b)}=0,
\end{equation}
where the second and third equalities were shown in the previous paragraph. 

We wish to isolate the
$\lambda_{j,w(j^\prime)}z_{b,w(j^\prime)}$ term on the lefthand side of (\ref{eqn:zerosumeqn1}) by showing all
other terms there vanish.  

\begin{Subclaim}
\label{subclaim:ffqq2}
For all $a \neq w(j^\prime)$ on the lefthand side
of (\ref{eqn:zerosumeqn1}), either $\lambda_{j,a} = 0$ or
$z_{b,a} = 0$.
\end{Subclaim}
\noindent
\emph{Proof of Subclaim~\ref{subclaim:ffqq2}:}
Given $p+1\leq a\leq n$, we prove the Subclaim by case analysis on the
position of the pivot in column $a$.

Since $\beta_b=+$ and $p+1\leq a\leq n$, the
pivot of $M$ in column $a$ cannot be in row $b$. This leaves two other cases.

First, if the pivot of $M$ in column $a$ is in a row strictly north of
row $b$,  then this pivot occurs in row $\caI(w^{-1}(a))<b$.  Note that we cannot have $\caI(w^{-1}(a)) = j^\prime$, since this would imply that $a=w(j)$, contradicting the fact that $p+1 \leq a \leq n$ while $1 \leq w(j) \leq p$.  Furthermore, we cannot have $\caI(w^{-1}(a)) = j$, since this implies that $a = w(j^\prime)$, and we are excluding this case from consideration.  Thus we need only consider the cases where $\caI(w^{-1}(a)) < j$, where $j < \caI(w^{-1}(a) < j^\prime$, and where $j^\prime < \caI(w^{-1}(a)) < b$.  In the first case, we have that $\lambda_{j,a}=\lambda_{j,w(\caI(\caI(w^{-1}(a))))}=0$
by Claim \ref{claim:isoproof1}.  In the second case, we have $\lambda_{j,a} = 0$ by Claim \ref{claim:isoproof2}.  
And in the third case, we have that $\lambda_{j,a} = 0$ by the inductive hypothesis.

The second possibility is that the pivot of $M$ in column $a$ is in a row strictly south of row $b$.  But in this case (Z.6) implies that $z_{b,a}=0$, so we are done.

This completes the proof of Subclaim~\ref{subclaim:ffqq2}.
\qed

By Claim~\ref{claim:lambdakk}, we have
$\lambda_{j,w(j^\prime)}=-1/2$. Hence, by
Subclaim~\ref{subclaim:ffqq2},
(\ref{eqn:zerosumeqn1}) reduces to
\[\lambda_{j,w(j^\prime)}z_{b,w(j^\prime)}=-1/2 z_{b,w(j^\prime)}=0.\] 
Therefore, $z_{b,w(j^\prime)}=0$, as needed.

($\beta$ has a $-$ at $b$):  The reasoning here is very similar to the previous case.  We omit the details.

($\beta$ has the left end of a matching at $b$):
By (O.2) in row $b$ of $M$, there are $1$'s in positions $w(b)$ and $w(\mathcal{I}(b))$. In addition, by (Z.1) and
(Z.5), $M$ has $0$'s elsewhere in row $b$.
Thus we obtain the claim's assertion that
\[z_{b,w(j)}=z_{b,w(j^\prime)}=0.\] 
Using this, taking $k=j$ in \eqref{eqn:oureqn1} and examining row $b$, we see that 
\begin{equation}
\label{eqn:ffqq3}
\lambda_{j,w(b)} = \lambda_{j,w(\mathcal{I}(b))}.
\end{equation} 
Claim~\ref{claim:zbw}'s hypothesis that $b > j^\prime$, combined with
Claim~\ref{claim:pibasis}(II), indicates that $\lambda_{j,w(b)} = 0$. Hence,
by (\ref{eqn:ffqq3}) we obtain
the required $\lambda_{j,w(b)} = \lambda_{j,w(\mathcal{I}(b))} = 0$.

($\beta$ has the right end of a matching at $b$ and its left end 
is at ${\mathcal I}(b)<j$):  
Since this subcase assumes ${\mathcal I}(b)<j$, by
Claim~\ref{claim:isoproof1}, we have that 
\[\lambda_{j,w({\mathcal I}(b))}=\lambda_{j,w(b)}=0,\]
as desired.

By definition of $w$,
\begin{equation}
\label{eqn:wsconsteqn}
w({\mathcal I}(b))<w(j)\leq p \mbox{\ and $p+1\leq w(b)<w(j^\prime)$}.
\end{equation}
Now, $z_{b,w({\mathcal I}(b))}=-1$ by 
(O.3), and column $w(j)$ is right of column $w({\mathcal I}(b))$ by the first inequality of \eqref{eqn:wsconsteqn},
so by (Z.4), we know that $z_{b,w(j)}=0$.  Similarly, $z_{b,w(b)} = 1$ by (O.3), and
using the last inequality of \eqref{eqn:wsconsteqn} combined with (Z.8), we
see that $z_{b,w(j^\prime)}=0$ as well.

($\beta$ has the right end of a matching at $b$ and its left end $\mathcal{I}(b)$ satisfies $j<{\mathcal I}(b)<j^\prime$):
By Claim~\ref{claim:isoproof2} we know $\lambda_{j,w(b)}=\lambda_{j,w(\mathcal{I}(b))} = 0$.  Here,
the matchings $(\mathcal{I}(b)<b)$ and $(j<j^\prime)$ are in a $1212$ pattern, so by (Z.9),
\[z_{b,w(j)}=-z_{b,w(j^\prime)}.\]  
Thus we only need to see that $z_{b,w(j)}=0$.  For this, consider row $b$ of
\eqref{eqn:oureqn2}, where we take $k=j$.  Since $b > j^\prime$, using Claim~\ref{claim:pibasis}(II), (O.3), and (Z.4), we see that
\begin{equation}
\label{eqn:ffeee1}
 \left(\displaystyle\sum_{a=1}^{w(\mathcal{I}(b))-1} \lambda_{j,a} z_{b,a}\right) - \lambda_{j,w(\mathcal{I}(b))} = 0. 
\end{equation}
It was argued already, in the proof of Claim \ref{claim:isoproof2} (see (R.3.2.1) - (R.3.2.3)), that the left hand side of (\ref{eqn:ffeee1}) reduces
to $1/2z_{b,w(j)}-\lambda_{j,w(\caI(b))}$, and we have noted above that $\lambda_{j,w(\caI(b))} = 0$.  So we have $1/2z_{b,w(j)} = 0$, whence $z_{b,w(j)}=0$, as desired.

($\gb$ has the right end of a matching at $b$ and its left end $\mathcal{I}(b)$ satisfies $j^\prime < \mathcal{I}(b)$):
Consider row $\mathcal{I}(b)$ of \eqref{eqn:oureqn1}, where we take $k=j$. Since we assume
$j^\prime < \mathcal{I}(b)$, by (O.2), (Z.1) and (Z.5),
\eqref{eqn:oureqn1} reads
\begin{equation}
\label{eqn:qwww1}
\lambda_{j,w(b)} = \lambda_{j,w(\mathcal{I}(b))}.
\end{equation}
By Claim~\ref{claim:pibasis}(II) and the assumption
$j^\prime<\caI(b)$, we see that 
$\lambda_{j,w(\caI(b))}=0$, and so by (\ref{eqn:qwww1}),
\begin{equation}
\label{eqn:trqp1} 	
\lambda_{j,w(b)}=0.
\end{equation}

Let us show $z_{b,w(j)} = 0$.  For this, we first analyze
row $b$ of \eqref{eqn:oureqn2}, taking $k=j$. By (O.3), (Z.8) and
Claim~\ref{claim:pibasis}(II), this row is
\[ \displaystyle\sum_{a=1}^{w(\mathcal{I}(b))-1} \lambda_{j,a} z_{b,a} - \lambda_{j,w(\mathcal{I}(b))} = 0. \]
Since by \eqref{eqn:qwww1} and \eqref{eqn:trqp1} we know
$\lambda_{j,w(\caI(b))} = 0$, this equation reduces to
\begin{equation}
\label{eqn:trqp2}
 \displaystyle\sum_{a=1}^{w(\mathcal{I}(b))-1} \lambda_{j,a} z_{b,a} = 0. 
\end{equation}
By Claim~\ref{claim:repeatedly} (where $k=\caI(b)$), all $a$ indexing 
the summation in (\ref{eqn:trqp2}) satisfy
\[w^{-1}(a) < \mathcal{I}(b) < b.\] 
Therefore, all $\lambda_{j,a}$ except for $\lambda_{j,w(j)}$ vanish, either by Claim \ref{claim:isoproof1}, Claim \ref{claim:isoproof2},
or the inductive hypothesis.  In view of 
Claim~\ref{claim:lambdakk},
\[\lambda_{j,w(j)} z_{b,w(j)} = 1/2 z_{b,w(j)}\] 
\textit{does} appear in the lefthand side of
(\ref{eqn:trqp2}). Thus (\ref{eqn:trqp2}) actually states
\[1/2z_{b,w(j)} = 0,\] 
whence $z_{b,w(j)} = 0$.

We now show $z_{b,w(j^\prime)}=0$.  For this, consider row $b$ of
\eqref{eqn:oureqn1} (taking $k=j$), which says that
\[ 0 + \displaystyle\sum_{a=p+1}^{w(b)-1} \lambda_{j,a} z_{b,a} + \lambda_{j,w(b)} = \displaystyle\sum_{a=1}^{w(\mathcal{I}(b)) - 1} \lambda_{j,a} z_{b,a} - \lambda_{j,w(\mathcal{I}(b))}, \]
or equivalently, by  (\ref{eqn:qwww1})
and (\ref{eqn:trqp1}), 
\[ \displaystyle\sum_{a=p+1}^{w(b)-1} \lambda_{j,a} z_{b,a} = \displaystyle\sum_{a=1}^{w(\mathcal{I}(b)) - 1} \lambda_{j,a} z_{b,a}. \]
By \eqref{eqn:trqp2}, the right hand side of this equation is zero, so we are reduced to
\begin{equation}
\label{eqn:qqqm1}
 \displaystyle\sum_{a=p+1}^{w(b)-1} \lambda_{j,a} z_{b,a} = 0. 
\end{equation}
By Claim~\ref{claim:repeatedly}, all $a$ indexing the sum 
satisfy $\mathcal{I}(w^{-1}(a)) < b$. Thus, 
all $\lambda_{j,a}$ vanish with the exception of $\lambda_{j,w(j^\prime)}$, either by Claim \ref{claim:isoproof1}, Claim
\ref{claim:isoproof2}, or the inductive hypothesis.  Thus \eqref{eqn:qqqm1} reduces to 
\[\lambda_{j,w(j^\prime)} z_{b,w(j^\prime)} =
-1/2 z_{b,w(j^\prime)} = 0,\] 
by Claim~\ref{claim:lambdakk}.  Hence $z_{b,w(j^\prime)} = 0$.

This concludes the proof of Claim \ref{claim:zbw}.
\qed

Now, we turn to rows $j$ and $j^\prime$.  By (Z.1) and (Z.5), every entry in row $j$ other than $z_{j,w(j)} = z_{j,w(j^\prime)} = 1$ is $0$.  Furthermore, by (Z.4) and (Z.8), 
\[ z_{j^\prime,a} = 0 \text{\ \ except when $1 \leq a \leq w(j)$ or $p+1 \leq a \leq w(j^\prime)$}. \]
Thus what remains to complete the proof of Proposition \ref{prop:projection-injective} is the following:

\begin{Claim}[Case C:  Row zeroness]
\label{claim:zja}
$z_{j^\prime,a}=0$ for all $a$ satisfying either $1\leq a<w(j)$ or $p+1\leq a<w(j^\prime)$.
\end{Claim}
\noindent
\emph{Proof of Claim~\ref{claim:zja}:}
By Claim~\ref{claim:repeatedly}, it suffices to prove 
\[z_{j^\prime,w(b)}=0 \text{\ and $z_{j^\prime, w(\caI(b))}=0$ for all $b<j$.}\] 
 We prove this by \textit{reverse} induction on $b$.  Our inductive hypothesis is that 
\[z_{j^\prime,w(c)}=z_{j^\prime,w(\caI(c))}=0 \text{\ for all $c$ with $b<c<j$.}\]  

Now, in view of (O.3), (Z.1), and (Z.5), row $j$ of equation~\eqref{eqn:oureqn1} (taking $k=b$) gives 
\begin{equation}
\label{eqn:uza1}
\lambda_{b,w(j)}=\lambda_{b,w(j^\prime)}.
\end{equation}

Now consider row $j^\prime$ of \eqref{eqn:oureqn1}, again taking $k=b$.  
By (O.3), (Z.4) and (Z.8),
\begin{equation}
\label{eqn:left-right}
 \sum_{a=1}^{w(j)-1}\lambda_{b,a}z_{j^\prime,a} - \lambda_{b,w(j)} = \sum_{a=p+1}^{w(j^\prime)-1}\lambda_{b,a}z_{j^\prime,a} + \lambda_{b,w(j^\prime)}.
\end{equation}
Hence, by Claim~\ref{claim:pibasisminusone}(II), we have
\begin{equation}
\label{eqn:rowjprimezerofirst}
\sum_{a=1}^{w(j)-1}\lambda_{b,a}z_{j^\prime,a} - \lambda_{b,w(j)} = 0,
\end{equation}
so by \eqref{eqn:left-right},
\begin{equation}
\label{eqn:rowjprimezerosecond}
\sum_{a=p+1}^{w(j^\prime)-1}\lambda_{b,a}z_{j^\prime,a} + \lambda_{b,w(j^\prime)} = 0
\end{equation}
as well.

We now consider multiple cases depending on the value of $\gb_b$.

($\beta_b=+$): By Claim~\ref{claim:repeatedly}, for $a$ indexing the sum in \eqref{eqn:rowjprimezerosecond}, we have $\caI(w^{-1}(a)) < j$.
Now, if $\caI(w^{-1}(a)) < b$, then $\lambda_{b,a}=0$ by Claim~\ref{claim:isoproof1}.  If $b<\caI(w^{-1}(a))<j$, then $z_{j^\prime,a}=0$ by 
the inductive hypothesis.  Note that $\caI(w^{-1}(a))=b$ is not possible, since this would imply that $a=w(b)$, contradicting the fact that
$p+1 \leq a < w(j^\prime)$ while $1 \leq w(b) \leq p$.  Thus \eqref{eqn:rowjprimezerosecond} reduces to $\lambda_{b,w(j^\prime)}=0$. Then by
(\ref{eqn:uza1}), $\lambda_{b,w(j)}=0$ also. Then Equation~(\ref{eqn:rowjprimezerofirst}) gives 
\begin{equation}
\label{eqn:nov17abc123}
\sum_{a=1}^{w(j)-1} \lambda_{b,a}z_{j^\prime,a}=0.
\end{equation}
  
Similarly to the previous paragraph, 
\[\lambda_{b,a}=0 \text{\ if $w^{-1}(a)<b$, and $z_{j^\prime,a}=0$ if $w^{-1}(a)>b$.}\]  
Hence \eqref{eqn:nov17abc123} reduces to
$\lambda_{b,w(b)}z_{j^\prime,w(b)}=0$.  Since $\lambda_{b,w(b)}\neq0$ by Claim~\ref{claim:lambdakk+-}, we have $z_{j^\prime,w(b)}=0$.  Since $\gb_b=+$, $\caI(b)=b$, so $z_{j^\prime,w(\caI(b))}=0$, as desired.

($\beta_b=-$): This is similar to the previous case, except the roles of Equations~(\ref{eqn:rowjprimezerofirst}) and~(\ref{eqn:rowjprimezerosecond}) are switched.

($\beta_b$ is a right endpoint):  First, consider \eqref{eqn:rowjprimezerofirst}.  Since $1 \leq a < w(j)$, by Claim
\ref{claim:repeatedly} we have that $w^{-1}(a) < j$.  So consider the cases $w^{-1}(a) < \caI(b)$, $w^{-1}(a) = \caI(b)$, 
$\caI(b) < w^{-1}(a) < b$, and $b < w^{-1}(a) < j$.  (Note that $w^{-1}(a)=b$ is not possible, since $1 \leq a < w(j) \leq p$,
while $p+1 \leq w(b) \leq n$.)  In the first case, we have that $\lambda_{b,a}=0$ by Claim~\ref{claim:isoproof1}.  In the
second case, we simply have that $a=w(\caI(b))$.  In the third case, we have that $\lambda_{b,a}=0$ by Claim~\ref{claim:isoproof2}.
Finally, in the fourth case, we have that $z_{j^\prime,a}=0$ by the inductive hypothesis.  Thus by Claim~\ref{claim:lambdakk},
\eqref{eqn:rowjprimezerofirst} reduces to
\begin{equation}
\label{eqn:rowjprimezerofirstreduced}
-1/2z_{j^\prime,w(\caI(b))}-\lambda_{b,w(j)}=0.
\end{equation}

Now, consider \eqref{eqn:rowjprimezerosecond}.  Since $p+1 \leq a < w(j^\prime)$, Claim \ref{claim:repeatedly} says that
$\caI(w^{-1}(a)) < j$.  So consider the cases $\caI(w^{-1}(a)) < \caI(b)$, $\caI(w^{-1}(a))=\caI(b)$,
$\caI(b) < \caI(w^{-1}(a)) < b$, and $b < \caI(w^{-1}(a)) < j$.  (As above, $\caI(w^{-1}(a)) = b$ cannot occur, since this 
would imply that $a=w(\caI(b))$, but $p+1 \leq a < w(j^\prime)$ while $1 \leq w(\caI(b)) < w(j) \leq p$.)  In the first case,
$\lambda_{b,a}=0$ by Claim \ref{claim:isoproof1}.  In the second case, we simply have that $a=w(b)$.  In the third case,
$\lambda_{b,a}=0$ by Claim \ref{claim:isoproof2}.  Lastly, in the fourth case, we have that $z_{j^\prime,a} = 0$ by induction.  Thus
using Claim \ref{claim:lambdakk} again, \eqref{eqn:rowjprimezerosecond} reduces to
\begin{equation}
\label{eqn:rowjprimezerosecondreduced}
-1/2z_{j^\prime,w(b)}+\lambda_{b,w(j^\prime)}=0.
\end{equation}

We can make arguments similar to the preceding ones when examining rows $j$ and $j^\prime$ of \eqref{eqn:oureqn1}, but taking  $k=\caI(b)$ instead.  Doing so, analogously to \eqref{eqn:uza1}, we obtain
\begin{equation}
\label{eqn:uza2}
\lambda_{\caI(b),w(j)} = \lambda_{\caI(b),w(j^\prime)}.
\end{equation}

Analogously to \eqref{eqn:rowjprimezerofirstreduced} and \eqref{eqn:rowjprimezerosecondreduced}, we obtain
\begin{equation}
\label{eqn:rowbmatchjprimezerofirstreduced}
1/2z_{j^\prime,w(\caI(b))}-\lambda_{\caI(b),w(j)}=0
\end{equation}
and
\begin{equation}
\label{eqn:rowbmatchjprimezerosecondreduced}
-1/2z_{j^\prime,w(b)}+\lambda_{\caI(b),w(j^\prime)}=0.
\end{equation}

The only simultaneous solution to \eqref{eqn:uza1}, \eqref{eqn:uza2}, \eqref{eqn:rowjprimezerofirstreduced}, \eqref{eqn:rowjprimezerosecondreduced}, \eqref{eqn:rowbmatchjprimezerofirstreduced}, and \eqref{eqn:rowbmatchjprimezerosecondreduced} is
$$z_{j^\prime,w(\caI(b))}=z_{j^\prime,w(b)}=\lambda_{\caI(b),w(j)}=\lambda_{\caI(b),w(j^\prime)}=\lambda_{b,w(j)}=
\lambda_{b,w(j^\prime)}=0.$$

($\beta_b$ is a left endpoint, and its right end satisfies $j^\prime<\caI(b)$): We have $z_{j^\prime,w(b)}=z_{j^\prime,w(\caI(b))}=0$ by (Z.3) and (Z.7) applied to columns $w(b)$ and $w(\caI(b))$ respectively.

($\beta_b$ is a left endpoint, and its right end satisfies $j<\caI(b)<j^\prime$): Since the matchings $(b<\caI(b))$ and $(j<j^\prime)$ are in a $1212$ pattern, by (Z.9),
\begin{equation}
\label{eqn:1212patternz}
z_{j^\prime,w(b)}=-z_{j^\prime,w(\caI(b))}.
\end{equation}

By precisely the arguments that preceded \eqref{eqn:rowjprimezerofirstreduced}
and
\eqref{eqn:rowjprimezerosecondreduced}, we see that \eqref{eqn:rowjprimezerofirst} reduces to
\[ \lambda_{b, w(b)}z_{j^\prime,w(b)} - \lambda_{b, w(j)} = 0, \]
while \eqref{eqn:rowjprimezerosecond} reduces to 
\[ \lambda_{b,w(\caI(b))}z_{j^\prime, w(\caI(b))} + \lambda_{b,w(j^\prime)}=0.\]

Using Claim~\ref{claim:lambdakk}, \eqref{eqn:uza1}, and \eqref{eqn:1212patternz}, we have
\[1/2z_{j^\prime,w(b)}-\lambda_{b,w(j)}=0 \text{\  and $1/2z_{j^\prime,w(b)}+\lambda_{b,w(j)}=0$.}\]  
This implies $z_{j^\prime,w(b)}=\lambda_{b,w(j)}=0$. Then we also have $z_{j^\prime,w(\caI(b))}=0$ by \eqref{eqn:1212patternz}, as required.

($\beta_b$ is a left endpoint, and its right end satisfies $\caI(b)<j$): In this case, the necessary claims are already proved by induction, since $b<\caI(b)<j$.

This completes the proof of Claim~\ref{claim:zja}, and hence the proof of Proposition~\ref{prop:projection-injective}.
\qed

\section{Conjectures and problems}

\begin{Conjecture}\label{conj:radical}
 $\caI_{\gamma,\ga}$ is a radical ideal.
\end{Conjecture}

Conjecture \ref{conj:radical} has been verified using Macaulay 2 through $p+q=6$. This conjecture would follow from a solution to:

\begin{Problem}
Find a Gr\"{o}bner basis for $\caI_{\gamma,\ga}$ with
square-free lead terms.
\end{Problem}

The analogous problem for \emph{Kazhdan-Lusztig ideals}
was solved in \cite{Woo-Yong-Grobner-12}; it was shown that the defining generators of the ideal form a Gr\"{o}bner basis. In contrast, we have:

\begin{Example}
Let $\alpha=-1221+$ and $\gamma=12-+12$. The defining generators of $\caI_{\gamma,\ga}$ are
\[z_{4,4}z_{6,6}+z_{5,4}z_{6,5}-z_{6,4}, \
-\frac{1}{8}z_{4,4}z_{6,6}-\frac{1}{8}z_{5,4}z_{6,5}+\frac{1}{4}z_{6,4}.\]
Thus $z_{6,4}\in \caI_{\gamma,\ga}$, and hence
$z_{4,4}z_{6,4}+z_{5,4}z_{6,5}\in \caI_{\gamma,\ga}$. Under any term order, the initial ideal must contain
$z_{6,4}$ and either $z_{4,4}z_{6,4}$ or $z_{5,4}z_{6,5}$. However, at most one of these monomials can be realized as a
lead term of the defining generators (for a fixed choice of term order). Hence the defining generators cannot be
a Gr\"{o}bner basis under any term order.\qed
\end{Example}

	The {\bf maximal singular locus} ${\tt Maxsing}(Y_{\gamma})$ is the set of $K$-orbits ${\mathcal O}_{\alpha}$
such that $Y_{\gamma}$ is singular along $\caO_{\ga}$ and $\caO_{\ga}$ is maximal in Bruhat order with respect to this property.

Whereas the theorem of W.~M.~McGovern \cite{McGovern-09a} recalled in the Introduction combinatorially characterizes
which $Y_{\gamma}$ are singular, the following question is open:
\begin{Problem}
\label{prob:singlocus}
Give a combinatorial description of ${\tt Maxsing}(Y_{\gamma})$.
\end{Problem}

	In \cite{Woo-Wyser-14}, a solution to Problem~\ref{prob:singlocus} is given for clans $\gamma$ which are
	``$1212$-avoiding''.  This gives a complete solution for the case $K=GL_{n-1} \times GL_1$
	(since all $(n-1,1)$-clans are $1212$-avoiding).  For the case of $K=GL_{n-2} \times GL_2$, we have a conjectural
	solution, but it is lengthy to state, so we omit it here.

The orbit closure $Y_{\gamma}$ is {\bf rationally smooth} along
${\mathcal O}_{\alpha}$ if the Kazhdan-Lusztig-Vogan polynomial    $P_{\gamma,\alpha}(q)$ equals $1$, and 
$Y_{\gamma}$ is (globally) rationally smooth if 
$P_{\gamma,\alpha}(q)=1$ for every $\alpha\leq \gamma$. It is known \cite{McGovern-09a}
that $Y_{\gamma}$ is rationally smooth if and only if $Y_{\gamma}$ is smooth. Using data from the {\tt ATLAS} project, we have verified the following conjecture for
$(p,q)=(2,2),(3,2)$:

\begin{Conjecture}\label{conj:ratl-smooth-locus}
For $\ga \leq \gamma$, $Y_{\gamma}$ is rationally smooth on ${\mathcal O}_{\alpha}$
if and only if $Y_{\gamma}$ is smooth on ${\mathcal O}_{\alpha}$.
\end{Conjecture}

Let $(R,{\mathfrak m})$ be the local ring of a point $p$ in a projective
variety $X$. The {\bf associated graded ring} is
\[{\rm gr}_{\mathfrak m}R=\bigoplus_{i\geq 0}{\mathfrak m}^i/{\mathfrak m}^{i+1}.\]

\begin{Conjecture}
Let $(R,{\mathfrak m})$ be the local ring associated to any point
$p\in {\mathcal O}_{\alpha}\subseteq Y_{\gamma}$. Then 
${\rm gr}_{{\mathfrak m}}R$ is Cohen-Macaulay and reduced.
Moreover, it is Gorenstein
whenever $Y_{\gamma}$ is Gorenstein along ${\mathcal O}_{\alpha}$.
\end{Conjecture}
This has been checked for $(p,q)=(2,2)$, $(3,2),(3,3),(4,3)$.
This conjecture does not follow from Conjecture~\ref{conj:radical}, 
since
reducedness, Cohen-Macaulayness, and Gorensteinness may be lost on
degenerating to the associated graded ring.

The {\bf projectivized tangent cone is} ${\rm Proj}({\rm gr}_{\mathfrak m}R)$.
The {\bf (Hilbert-Samuel) multiplicity} of $p\in X$ is
\[{\rm mult}_p(X)={\rm deg}({\rm Proj}({\rm gr}_{\mathfrak m}R)).\]
This statistic provides singularity information. Specifically, ${\rm mult}_p(X)=1$ if and only if $p$ is a smooth point of $X$.
Let ${\rm mult}_{\gamma,\alpha}$ be
${\rm mult}_p(Y_{\gamma})$ for any point $p\in {\mathcal O}_{\alpha}$.

 Define 
\[{\tt maxmult}(Y_{\gamma})=\max_{\alpha\leq\gamma} {\tt mult}_{\gamma,\alpha}.\] 
Since multiplicity is a semicontinuous numerical invariant, 
${\tt maxmult}(Y_{\gamma})={\tt mult}_{\gamma,\alpha}$ for some
matchless $\alpha\leq\gamma$ (so ${\mathcal O}_{\alpha}$ is a closed
orbit). A search for a combinatorial rule for  ${\rm mult}_{\gamma,\alpha}$
might be partially guided by a solution to the following problem:

\begin{Problem}
\label{conj:maxmult}
Is the maximum value of ${\tt maxmult}(Y_{\gamma})$ for $\gamma\in {\tt Clans}_{p,q}$
achieved at $\gamma_{\rm max}=1+^{p-1}-^{q-1}1\in {\tt Clans}_{p,q}$? 
\end{Problem}

We have checked that the answer to Problem~\ref{conj:maxmult} is affirmative $p+q\leq 6$. 
Note the maximizer in Problem~\ref{conj:maxmult} 
need not be unique. For instance, when $p=q=2$, the maximum is achieved at
$\gamma=1212$, $\gamma=1+-1$, and $\gamma=1-+1$, with ${\tt maxmult}(Y_{\gamma})=2$.

\begin{figure}[h!]
        \centering
        \includegraphics[scale=0.7]{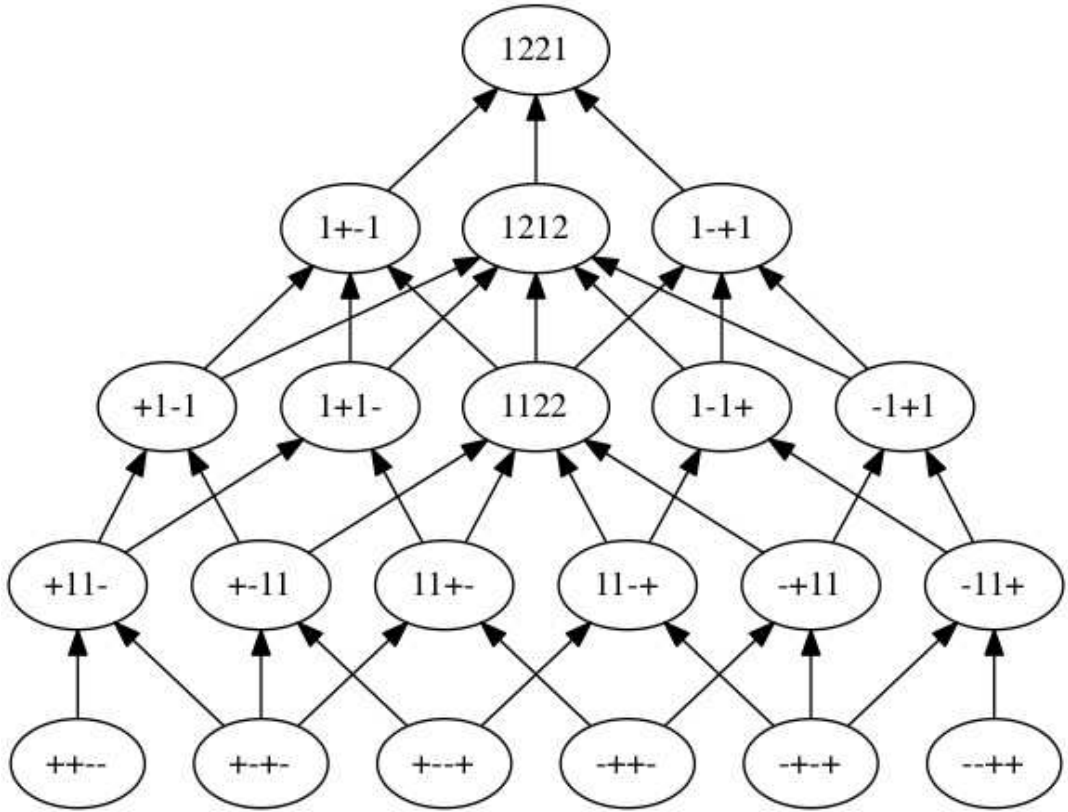}
        \caption{Bruhat order on $(GL(4,\C),GL(2,\C) \times
GL(2,\C))$}\label{fig:type-a-2-2}
\end{figure}

\begin{Example}[Analogue of G.~Lusztig's conjecture is false]
\label{exa:Lusztig}
Note that $P_{1+-1,++--}(q)=q+1$, whereas $P_{1212,++--}(q)=1$. However, 
$[++--,1+-1]\cong [++--,1212]$ as posets (see Figure~\ref{fig:type-a-2-2}). Thus, the KLV polynomial is not an invariant of the poset in Bruhat order. This contrasts with
the situation for Schubert varieties, where a conjecture attributed to G.~Lusztig asserts that the
Kazhdan-Lusztig polynomial is a poset invariant; see for example \cite{Brenti-03} and the references
therein. Similarly, 
${\tt mult}_{1+-1,++--}=2$ whereas ${\tt mult}_{1212,++--}=1$.\qed
\end{Example}

The {\bf Cohen-Macaulay type} of a local Cohen-Macaulay ring $R$ is $\dim_{\Bbbk} {\rm Ext}_R^{\dim R}({\Bbbk},R)$.  A ring $R$ is
{\bf Gorenstein} if its Cohen-Macaulay type is $1$. 
\excise{In general if $R=S/I$
where $S$ is a polynomial ring and $I\subset S$ is an ideal, one can compute the type using the
following {\tt Macaulay 2} commands:
\begin{verbatim}
     S:=ring I;
     R:=S/I;
     k=coker vars R;
     numgens Ext^(dim(S)-dim(I))(S^1/I,S^1);
\end{verbatim}
(See, e.g., \textsf{http://mathoverflow.net/questions/122118/how-to-check-if-a-commutative-ring
-is-gorenstein}.)

Since $Y_{\gamma}$ is Cohen-Macaulay, it follows that if $I$ is
the ideal $\sqrt{I_{\gamma,\ga}}$, the
above commands will determine the Cohen-Macaulay type at any point of a given orbit ${\mathcal O}_{\alpha} \subseteq Y_{\gamma}$.}
	
	Define the {\bf maximal non-Gorenstein locus} 
${\tt MaxnonGor}(Y_{\gamma})=\{{\mathcal O}_{\alpha}\}$ to be the set
of $K$-orbits ${\mathcal O}_{\alpha} \subseteq Y_{\gamma}$ such that $Y_\gamma$ is non-Gorenstein at any point of $\mathcal{O}_{\ga}$
and $\caO_{\ga}$ is maximal in Bruhat order with respect to this property.

\begin{Conjecture}\label{conj:gorenstein-sing-locus}
${\tt MaxnonGor}(Y_{\gamma})\subseteq {\tt Maxsing}(Y_{\gamma})$.
\end{Conjecture}

If $Y_{\gamma}$ is non-Gorenstein along $\mathcal{O}_{\ga}\in
{\tt Maxsing}(Y_{\gamma})$ then $\mathcal{O}_{\ga}\in {\tt MaxnonGor}(Y_{\gamma})$.
However, if $Y_\gamma$ is Gorenstein along some
$\mathcal{O}_{\ga}\in {\tt Maxsing}(Y_{\gamma})$, there a priori may be 
$\mathcal{O}_{\beta}\in {\tt MaxnonGor}(Y_{\gamma})$ with $\beta < \ga$.  Conjecture \ref{conj:gorenstein-sing-locus} asserts 
that this does not occur.  This conjecture is an analogue of \cite[Conjecture 6.7]{Woo-Yong-08}.

We verified Conjecture \ref{conj:gorenstein-sing-locus} for all $(p,q)$ with $p+q\leq 7$ and for many cases of $(p,q)=(4,4)$. 
Below, we give the singularity data 
$(\gamma,\ell(\gamma), {\tt Maxsing}(Y_{\gamma}), {\tt MaxnonGor}(Y_{\gamma}))$
for $p=3, q=2$ for all the clans $\gamma$ where $Y_\gamma$ is singular.

\begin{verbatim}
(122+1, 5, {1-1++}, {})
(1+221, 5, {++1-1}, {})
(12+12, 5, {-+++-, +1-1+}, {})
(1+212, 4, {++--+, 11++-}, {++--+})
(1-++1, 4, {--+++}, {--+++})
(121+2, 4, {+--++, -++11}, {+--++})
(1+-+1, 4, {11-++, ++-11}, {})
(1++-1, 4, {+++--}, {+++--})
(1+-1+, 3, {++--+}, {})
(+1+-1, 3, {+++--}, {})
(+1-+1, 3, {+--++}, {})
(1212+, 3, {+--++, -++-+}, {})
(+1212, 3, {+-++-, ++--+}, {})
(1-+1+, 3, {--+++}, {})
\end{verbatim}

\begin{Problem}
When is $Y_{\gamma}$ (globally) Gorenstein?
\end{Problem}

An answer to this question has been given
by the first two authors in \cite{Woo-Wyser-14} in the event that $\gamma$ is $1212$-avoiding.  As mentioned in the
Introduction, it is shown in \emph{loc. cit.} that Gorensteinness cannot be characterized by ordinary pattern avoidance
in general.  However, for small $q$, it apparently can.

When $q=1$, all orbit closures are smooth
(hence Gorenstein).  Any $(p,1)$-clan ($p \geq 1$) has at most one matching, and if it is has a matching, it can
have no $-$ signs.  Thus it necessarily avoids all of McGovern's singular patterns recalled in the
Introduction.

Next consider the case $q=2$.  In this case, the Gorensteinness criterion of \cite{Woo-Wyser-14} is easily
seen to amount to the following ordinary pattern avoidance criterion.

\begin{Fact}[{\cite[Proposition 3.4.2]{Woo-Wyser-14}}]
 If $\gamma$ is a $1212$-avoiding $(p,2)$-clan with $p \geq 2$, then $Y_{\gamma}$ is
 Gorenstein if and only if $\gamma$ avoids $1++-1$, $1-++1$, $1++221$, and $122++1$.
\end{Fact}

Experimentally, it seems that when $q=2$, Gorensteinness can be characterized by ordinary pattern
avoidance even if we allow $1212$-including clans.

\begin{Conjecture}\label{conj:gorenstein-patterns}
If $\gamma$ is any $(p,2)$-clan with $p \geq 2$, then $Y_{\gamma}$ is Gorenstein if and only if
$\gamma$ avoids $1++-1$, $1-++1$, $1++221$, $122++1$, $1+212$, and $121+2$.
\end{Conjecture}

This conjecture holds (by computation) for $n \leq 7$.  The example recalled in the Introduction says that
there is no ordinary pattern avoidance characterization when $q > 2$.

Recall that a local ring $R$ is said to be a {\bf local complete intersection (lci)}
if it is the quotient of a regular local ring by an 
ideal generated by a regular sequence. A variety is lci if
each of its local rings are lci. Every smooth variety is lci, while every lci variety is Gorenstein (and hence
Cohen-Macaulay).

At present, we do not have a characterization of which $Y_{\gamma}$ are
lci in general.  An ordinary pattern avoidance criterion is given in \cite{Woo-Wyser-14} for $1212$-avoiding clans.  For
general $p$ and $q$, there are many non-lci patterns that must be avoided.  However, when $q = 2$, the
list of non-lci patterns that can actually occur is much smaller.

\begin{Fact}[{\cite[Proposition 3.3.3]{Woo-Wyser-14}}]
  If $\gamma$ is a $1212$-avoiding $(p,2)$-clan with $p \geq 2$, then $Y_{\gamma}$ is lci if and 
  only if $\gamma$ avoids $1++-1$, $1-++1$, $1++221$, and $122++1$.  Thus $Y_{\gamma}$ is lci if and only if it is
  Gorenstein.
\end{Fact}

It is conjectured in \cite{Woo-Wyser-14} that ordinary pattern avoidance can be used to characterize lci-ness in general,
whether $\gamma$ is $1212$-avoiding or not.  The complete list of non-lci patterns is not known, even conjecturally.
However, for $q=2$, we have:

\begin{Conjecture}
 If $\gamma$ is any $(p,2)$-clan with $p \geq 2$, then $Y_{\gamma}$ is lci if and only if $\gamma$ avoids
 $1++-1$, $1-++1$, $1++221$, $122++1$, $1+212$, and $121+2$.  Equivalently, we conjecture that $Y_{\gamma}$ is lci if and
 only if it is Gorenstein.
\end{Conjecture}

This conjecture has been checked for $n \leq 7$.  The equivalence of lci-ness and Gorensteinness breaks down when
$q > 2$; for example, $Y_{122331}$ is Gorenstein but not lci.

\section*{Acknowledgements}
We thank Frank Sottile for a conversation years ago about \cite{Woo-Yong-08}
that ultimately led to this project.
We thank Michel Brion for many helpful comments and in particular for his suggestion to consider the slices
of \cite{Mars-Springer}, as well as his
assistance with Lemmas~\ref{lem:slice-stable-properties} and~\ref{lem:slice-of-orbit-closure}. 
AW was supported by NSA grant H98230-13-1-0242 and Simons Collaboration Grant 359792.  BW was supported by NSF International Research Fellowship 1159045
and hosted by Institut Fourier in Grenoble. AY was supported by NSF grants DMS 1201595 and DMS 1500691. This work was partially completed while AY and AW were visitors at the Fields Institute program on ``Combinatorial Algebraic Geometry''; they thank the organizers and the Institute for their hospitality.

\bibliography{sourceDatabase}
\bibliographystyle{plain}

\end{document}